\documentclass[12pt,twoside]{amsart}
\usepackage{amssymb}
\usepackage{amscd}

\title[Higher direct images]{Higher direct images of 
log canonical divisors and positivity theorems}
\author{Osamu Fujino} 
\date{2003/2/7}
\address{Research Institute for Mathematical Sciences\\ 
 Kyoto University, Kyoto 606-8502 Japan}
\email{fujino@kurims.kyoto-u.ac.jp}
\newcommand{\bQ}[0]{{\mathbb Q}}
\newcommand{\tor}[0]{\operatorname{tor}}
\newcommand{\xdiv}[0]{{\operatorname{div}}}
\newcommand{\xdiscrep}[0]{{\operatorname{discrep}}}
\newcommand{\xlog}[0]{{\operatorname{log}}}
\newcommand{\Supp}[0]{{\operatorname{Supp}}}
\newcommand{\codim}[0]{{\operatorname{codim}}}
\newcommand{\Exc}[0]{{\operatorname{Exc}}}
\newcommand{\Gr}[0]{{\operatorname{Gr}}}
\newcommand{\Ker}[0]{{\operatorname{Ker}}}
\newcommand{\Coker}[0]{{\operatorname{Coker}}}
\newcommand{\xO}[0]{{\operatorname{\mathcal O}}}
\newcommand{\xN}[0]{{\operatorname{\mathbb Z_{>0}}}}
\newtheorem{thm}{Theorem}[subsection]
\newtheorem*{theorem}{Theorem}         %\renewcommand{\thetheorem}{} 
 
\newtheorem{lem}[thm]{Lemma}
\newtheorem{cor}[thm]{Corollary}
\newtheorem{prop}[thm]{Proposition}

\theoremstyle{definition}
\newtheorem{defn}[thm]{Definition}
\newtheorem{cau}[thm]{Caution}
\newtheorem{rem}[thm]{Remark}
\newtheorem*{ack}{Acknowledgments}         %\renewcommand{\theack}{} 

\newtheorem{ex}[thm]{Example}
\theoremstyle{remark}
\newtheorem*{claim}{\bf{Claim}}  
  
\newtheorem{say}[thm]{\bf{}}
\newtheorem{step}{\bf{Step}}

\newtheorem{thm5}{\bf{Theorem}}[section]
\newtheorem{lem5}[thm5]{\bf{Lemma}}
\newtheorem{prop5}[thm5]{\bf{Proposition}}

\theoremstyle{definition}
\newtheorem{defn5}[thm5]{Definition}
\newtheorem{rem5}[thm5]{Remark}
\theoremstyle{remark}
\newtheorem{say5}[thm5]{\bf{}}

\newtheorem{thma}{\bf{Theorem}}
\newtheorem{propa}[thma]{\bf{Proposition}}
\begin{document}
\bibliographystyle{amsalpha+}

\begin{abstract} 
In this paper, we investigate higher direct images of 
log canonical divisors. 
After we reformulate Koll\'ar's torsion-free theorem, 
we treat the relationship between 
higher direct images of log canonical divisors and 
the canonical extensions of Hodge filtration of 
gradedly polarized variations of mixed Hodge structures. 
As a corollary, we obtain a logarithmic version of 
Fujita-Kawamata's semi-positivity theorem. 
By this semi-positivity theorem, we 
generalize Kawamata's positivity theorem and 
apply it to the study of a log canonical bundle formula. 
The final section is an appendix, which is a result of Morihiko 
Saito. 
\end{abstract}

\maketitle
\tableofcontents

\section{Introduction}\label{sec1}

In this paper, we investigate higher direct images of 
log canonical divisors. 

First, we reformulate Koll\'ar's torsion-free theorem 
and vanishing theorem. 
This part is more or less known to experts. 
See \cite{ev} and \cite[Section 3]{am}. 
However, we explain the details since 
there are no appropriate references for our purposes 
and torsion-freeness will play important roles in this 
paper. 

Next, we treat the relationship between 
higher direct images of log canonical divisors and 
the canonical extensions of Hodge filtration of 
gradedly polarized variations of mixed Hodge structures. 

Let $f:X\longrightarrow Y$ be a surjective morphism 
between non-singular projective varieties and 
$D$ a simple normal crossing divisor on X. 
We assume that $D$ is {\em{strongly horizontal}} (see 
Definition \ref{strong}) with respect to 
$f$. Then, under some suitable assumptions, 
$R^if_*\omega_{X/Y}(D)$ is characterized 
as the (upper) canonical extension of the 
bottom Hodge filtration of the suitable 
polarized variation of mixed Hodge structures. 
When $D=0$, it is the theorem of Koll\'ar and Nakayama 
(see \cite[Theorem 2.6]{ko2} and \cite[Theorem 1]{n}). 
If $Y$ is a curve, then the above theorem immediately follows 
from the study of the gradedly polarized variation of 
mixed Hodge structures by Steenbrink and Zucker (see 
\cite[\S 5 The geometric case]{sz}). 
By this characterization, it is not difficult to see 
that $R^if_*\omega_{X/Y}(D)$ is semi-positive on some monodromy 
conditions. 
It is a logarithmic version of Fujita-Kawamata's semi-positivity 
theorem. 

Finally, by using this semi-positivity theorem, we 
generalize Kawamata's positivity theorem. 
This is one of the main purposes of this paper. 
As a corollary, we obtain a log canonical bundle 
formula for log canonical pairs, 
which is a slight generalization of 
\cite[Section 4]{fm}. 

We don't pursue further applications or generalizations 
to make this paper readable. 

The final section is an appendix, which is a result of 
Morihiko Saito. 

\subsection{Main Results}\label{subsec1.1}

Let us explain the results of this paper more precisely. 
We will work over $\mathbb C$, the complex number field, 
throughout this paper. 

\begin{say}\label{01}
In Section \ref{sec2}, we  reformulate Koll\'ar's torsion-free 
theorem. 

\begin{theorem}[{cf.~Theorems \ref{torsionfree}, \ref{hosoku}}] 
Let $f:X\longrightarrow Y$ be a surjective morphism between 
projective varieties. Assume that 
$X$ is non-singular and $D:=\sum _{i\in I}D_i$ a 
simple normal crossing divisor on $X$. 
We assume that $D$ is {\em{strongly horizontal}} 
with respect to $f$, that is, 
every irreducible component of $D_{i_1}\cap D_{i_2}\cap 
\cdots \cap D_{i_k}$, where $\{i_1, \cdots, i_k\}\subset I$, 
is dominant onto $Y$ $($see {\em{Definition \ref{strong}}}$)$. 
Then $R^if_*\omega_{X/Y}(D)$ is torsion-free. 
\end{theorem}

It is a special case of \cite[Theorem 3.2 (i)]{am}. 
His theorem is much more general than ours. 
We explain the details and give a precise proof. 
Our proof is a modification of Arapura's 
argument \cite[Theorem 1]{a} and 
relies on the theory of (geometric) variation of 
mixed Hodge structures over curves. 
So, it is a warm-up to the next section. 

In subsection \ref{subsec2.2}, we treat a slight 
generalization of Koll\'ar's vanishing 
theorem (see Theorem \ref{vanishing}). 
We note that we don't use it later. 
Thus, we omit it here. 
\end{say}

\begin{say}\label{02}
Section \ref{sec3} is one of the main parts of 
this paper. It is a logarithmic generalization 
of the theorem of Koll\'ar and Nakayama. 
As a corollary, we obtain a logarithmic generalization of 
Fujita-Kawamata's semi-positivity theorem. 

\begin{theorem}[{cf.~Theorems \ref{can-ext}, \ref{loc-free}}] 
Let $f:X\longrightarrow Y$ be a surjective morphism 
between non-singular projective varieties and 
$D$ a simple normal crossing divisor 
on $X$, which is strongly horizontal with 
respect to $f$. 
Let $\Sigma$ be a simple normal crossing divisor on $Y$. 
We put $Y_0:=Y\setminus \Sigma$. 
If $f$ is smooth and $D$ is relatively normal crossing 
over $Y_0$, then $R^if_*\omega_{X/Y}(D)$ is the 
upper canonical extension of the bottom Hodge 
filtration. In particular, it is locally free. 
\end{theorem}

Note that on the above assumptions we 
have a (geometric) variation of mixed Hodge structures 
on $Y_0$. 
Our theorem is a direct consequence of 
\cite[\S 5]{sz} when $Y$ is a curve. 
If $D=0$, then it is the theorem of Koll\'ar and Nakayama 
(see \cite[Theorem 2.6]{ko2} and \cite[Theorem 1]{n}). 
A key point of our proof is the torsion-freeness of 
$R^if_*\omega_{X/Y}(D)$ that is obtained in Section \ref{sec2}. 

We put $X_0:=f^{-1}(Y_0)$, $D_0:=D\cap X_0$, 
$f_0:=f|_{X_0}$, and $d:=\dim X-\dim Y$. 

\begin{theorem}[{cf.~Theorem \ref{semi}}] 
We further assume that all the local monodromies on the 
local system $R^{d+i}
{f_0}_*\mathbb C_{X_0-D_0}$ around every irreducible 
component of $\Sigma$ are unipotent, 
then $R^if_*\omega_{X/Y}(D)$ is a semi-positive 
vector bundle. 
\end{theorem}

As stated above, 
it is a logarithmic version of Fujita-Kawamata's semi-positivity 
theorem. 
This theorem will play crucial roles in Section \ref{sec4}. 
\end{say}

\begin{say}\label{03}
Section \ref{sec4} deals with 
a generalization of Kawamata's positivity theorem \cite[Theorem 2]{ka3}. 
The statement is too technical to state here. 
So, please see Theorem \ref{kpt}. 
Our proof is essentially the same as Kawamata's. 
In his proof, 
he used Fujita-Kawamata's semi-positivity theorem. 
On the other hand, we use the semi-positivity of $f_*\omega_{X/Y}(D)$, 
which is obtained in Section \ref{sec3}. 
Roughly speaking, Kawamata's original positivity 
theorem holds for (sub) klt (Kawamata log terminal) pairs 
and our theorem (Theorem \ref{kpt}) does for 
(sub) lc (log canonical) pairs. 
We believe that this difference is big for some applications. 
We note that Kawamata's original positivity theorem 
already played important roles in 
various situations. 

In subsection \ref{subsec4.2}, 
we treat only one application of our positivity 
theorem. It is a slight generalization of 
\cite[Theorem 0.2]{f1}. 
See Theorem \ref{mt}. 
We don't pursue other applications here. 
\end{say}

\begin{say}\label{04}
In Section \ref{sec5}, we formulate and prove a log canonical 
bundle formula for lc pairs. 
For the precise formula, see Theorem \ref{mainthm}. 
This section is essentially the same 
as \cite[Section 4]{fm}, where we formulate and 
prove it 
for klt pairs. 
The only one nontrivial point is the semi-positivity (nefness) 
of the {\em{log-semistable part}} $L^{\log,ss}_{X/Y}$ 
(see Theorem \ref{nef}). 
It is a direct consequence of the positivity theorem obtained 
in Section \ref{sec4}. 

When we wrote \cite[Section 4]{fm}, Kawamata's positivity theorem 
was proved only for 
(sub) klt pairs. 
So we formulated a log canonical bundle formula for klt pairs. 
Since we generalized Kawamata's positivity 
theorem in Section \ref{sec4}, there are no difficulties to formulate 
and prove a log canonical bundle formula for lc pairs. 
We repeat the formulation in details for the readers' convenience. 
We note that it is conjectured that 
the log-semistable part $L^{\log, ss}_{X/Y}$ is semi-ample. 
It is proved only for elliptic fibrations, Abelian fibrations, 
$K3$ fibrations and so on. 
We recommend the readers to see \cite[Section 6]{f2} for 
the details. 

Note that F.~Ambro treated a log canonical bundle formula 
in a slightly different formulation. 
We don't pursue this formulation in this paper. 
For the details, see his preprint \cite{am2}. 
Related topics are \cite{am} and \cite{sho}.  
\end{say}

\begin{say}\label{ap}
Section \ref{sec6} is an appendix, which is a remark on 
Section \ref{sec3}. 
After I finished the preliminary version of 
this paper, I asked Professor Morihiko Saito 
about the topic in Section \ref{sec3}. 
I received an e-mail \cite{e-s} from him, 
where he gave a different proof (Proposition \ref{p2} in \ref{6.1}) to 
Theorems \ref{can-ext} and \ref{loc-free}. 
It depends on the theory of mixed Hodge Modules \cite{mhm}, 
\cite{kc}. 
I insert it into this paper as an appendix. 
Note that I made no contribution to Section \ref{sec6}. 
\end{say}

\begin{say}\label{00}
Subsection \ref{subsec1} collects some basic definitions 
and fix our notation. We also 
recall some vanishing theorems. 
After checking subsection \ref{subsec1} quickly, 
the readers can 
read any section independently 
with referring results obtained in other sections. 
\end{say}

\begin{ack}
I was partially supported by the Inoue Foundation for Science. 
I would like to express my sincere gratitude to Professors 
Yoshio Fujimoto, Masaki Kashiwara, 
Noboru Nakayama, and Hiromichi Takagi for 
useful comments. 
I am grateful to Professor Morihiko Saito very much, 
who kindly told 
me his result \cite{e-s} (see \ref{ap} above) 
and allowed me to use it. 
He also pointed out some ambiguities in the 
proof of Theorem \ref{can-ext}. 
His comments made this paper more readable. 
I am also grateful to Professor Shigefumi Mori 
for warm encouragements. 
\end{ack}
 
\subsection{Preliminaries}\label{subsec1} 

Let us recall the basic definitions and fix our 
notation (cf.~\cite{KMM}, \cite{kom}, and \cite{FA}). 
We also recall some vanishing theorems. 
We will work over $\mathbb C$, the complex number field, throughout 
this paper. 

\begin{defn}[{$\mathbb Q$-divisors}]\label{div} 
Let $X$ be a normal variety and $B, B'$ $\bQ$-divisors on $X$. 
If $B-B'$ is effective, we write $B \succ B'$ or $B'\prec B$. 
We write $B \sim B'$ if $B-B'$ is
a principal divisor on $X$ (linear equivalence of $\bQ$-divisors). 
Let $B_+, B_-$ be the effective $\bQ$-divisors on $X$
without common irreducible components such that 
$B_+-B_-=B$. They are called the {\it positive}
and the {\it negative} parts of $B$.
\end{defn}

\begin{defn}[Operations of $\mathbb Q$-divisors]\label{ope}
Let $B=\sum b_iB_i$ be a $\mathbb Q$-divisor, 
where $b_i$ are rational numbers and $B_i$ are 
mutually prime Weil divisors. 
We define 
\begin{eqnarray*}
\llcorner B\lrcorner &:=&\sum \llcorner b_i\lrcorner B_i, \ {\text{the 
round down of}} \ B, \\
\ulcorner B\urcorner &:=&\sum \ulcorner b_i\urcorner B_i=
-\llcorner -B\lrcorner, \ {\text{the 
round up of}} \ B, \\
B^{<1}&:=&\sum _{b_i<1}b_i B_i, 
\end{eqnarray*}
where for $r\in \mathbb R$, we define $\llcorner r \lrcorner :=
\max \{t\in \mathbb Z; t\leq r\}$. 
\end{defn}

\begin{defn}[Vertical and horizontal]\label{ver-hori}
Let $f : X \longrightarrow S$ be a surjective morphism between varieties.
Let $B^h, B^v$ be the $\bQ$-divisors on $X$
with $B^h+B^v=B$
such that an irreducible component of $\Supp B$
is contained in $\Supp B^h$ if and only if 
it is mapped onto $S$.
They are called the {\it horizontal}
and the {\it vertical} parts of $B$ over $S$.
A $\mathbb Q$-divisor 
$B$ is said to be {\it horizontal} (resp.~{\it vertical}\ )
over $S$ if $B=B^h$ (resp.~$B=B^v$). 
The phrase ``over $S$" might be suppressed
if there is no danger of confusion.
\end{defn}

\begin{defn}[Canonical divisor]\label{notation1}
Let $X$ be a normal variety. 
The {\em{canonical divisor}} $K_X$ is defined so that 
its restriction to the regular part of $X$ is a divisor of a regular 
$n$-form. 
The reflexive sheaf of rank one $\omega _X:={\mathcal O}_X(K_X)$ 
corresponding to $K_X$ is called the {\em{canonical sheaf}}. 
\end{defn}

The following is the definition of singularities of pairs. 
Note that the definitions in \cite {KMM} or \cite{kom} are 
slightly different from ours. 

\begin{defn}[Discrepancies and singularities for pairs]\label{def1}
Let $X$ be a normal variety and $D=\sum d_i D_i$ a 
$\bQ$-divisor on $X$ such that $K_X+D$ is $\bQ$-Cartier. 
Let $f:Y\longrightarrow X$ be a proper birational morphism 
from a normal variety $Y$. 
Then we can write 
$$
K_Y=f^{*}(K_X+D)+\sum a(E,X,D)E, 
$$ 
where the sum runs over all the distinct prime divisors $E\subset Y$, 
and $a(E,X,D)\in \bQ$. This $a(E,X,D)$ is called the 
{\em discrepancy} of $E$ with respect to $(X,D)$. 
We define 
$$
\xdiscrep (X,D):=\inf _{E}\{a(E,X,D)\ 
|\  E \text{ is exceptional over}\  X  \}.
$$
On the assumption that $d_i \leq 1$ 
for every $i$, we say that $(X,D)$ is 
$$
\begin{cases}
\text{sub klt}\\
\text{sub lc}\\
\end{cases}
\quad {\text{if}} \quad \xdiscrep (X,D) 
 \quad
\begin{cases}
>-1\quad {\text {and \quad $\llcorner D\lrcorner \leq 0$,}}\\
\geq -1.\\
\end{cases}
$$ 
If $(X,D)$ is sub klt (resp.~sub lc) and $D$ is effective, 
then we say that $(X,D)$ is {\em{klt}} (resp.~{\em{lc}}). 
Here klt (resp.~lc) is short for {\em {Kawamata log terminal}} (resp.~
{\em{log canonical}}). 
\end{defn}

\begin{defn}[Center of lc singularities]\label{def2}
Let $X$ be a normal variety and $D$ a $\mathbb Q$-divisor 
on $X$ such that $K_X+D$ is $\mathbb Q$-Cartier. 
A subvariety $W$ of $X$ is said to be a 
{\em center of log canonical singularities} 
for the pair $(X,D)$, if there exists a proper birational morphism 
from a normal variety $\mu:Y\longrightarrow X$ and a prime divisor $E$ on 
$Y$ with the discrepancy coefficient $a(E,X,D)\leq -1$ such 
that $\mu(E)=W$. 
\end{defn}

\begin{rem}[{cf.~\cite[Lemmas 2.29, 2.30, and 2.45]{kom}}]\label{compu} 
Let $X$ be a non-singular variety and $D$ a simple normal 
crossing divisor on $X$. 
Then $(X,D)$ is lc. More precisely, $(X,D)$ is a typical example of 
dlt pairs (see Remark \ref{dlt} below). 
Let $D=\sum _{i\in I}D_i$ be the irreducible decomposition of 
$D$. 
Then, $W$ is a center of log canonical singularities for the 
pair $(X,D)$ if and only if $W$ is an irreducible 
component of $D_{i_1}\cap D_{i_2}\cap \cdots \cap D_{i_k}$ 
for some $\{i_1, i_2, \cdots, i_k\}\subset I$. 
\end{rem}

\begin{defn}[Log canonical threshold]\label{lc-thre}
Let $(X,D)$ be a sub lc pair and 
$\Delta\ne 0$ a $\mathbb Q$-Cartier divisor on $X$. 
We put 
$$
c_0:=\max\{c\in \mathbb R\ | \ (X, D+c\Delta)\ {\text{is sub lc}}\}. 
$$ 
We call $c_0$ {\em{the log canonical threshold for the pair $(X,D)$ 
with respect to $\Delta$}}. 
We note that $c_0\in \mathbb Q$. 

Furthermore, we assume that $(X,D)$ is lc (resp.~klt) and 
$\Delta$ is an effective Weil divisor. 
Then $0\leq c_0\leq 1$ (resp.~$0<c_0\leq 1$). 
\end{defn}

We introduce the following new notion, which will play important 
roles in this paper. 

\begin{defn}[Strongly horizontal]\label{strong}
Let $(X, D)$ be a sub lc pair 
and $f:X\longrightarrow Y$ a surjective morphism. 
If all the center of log canonical singularities for the 
pair $(X,D)$ are dominant onto $Y$, then we call $D$ {\em{strongly 
horizontal with respect to $f$}}.  
\end{defn}

\begin{rem}\label{zero}
Let $f:X\longrightarrow Y$ and $D$ be as in Definition \ref{strong}. 
If $(X,D)$ is sub klt, then $D$ is strongly horizontal 
with respect to $f$. It is obvious by the definition. 
We note that $D=0$ is strongly horizontal when 
$X$ is non-singular.  
\end{rem}

We note the following easy fact. 

\begin{lem}\label{akiraka}
Let $f:X\longrightarrow Y$ and $D$ be as in 
{\em{Definition \ref{strong}}}, that is, 
$D$ is strongly horizontal with respect to $f$. 
Let $\Lambda$ be a free linear system on $Y$ and 
$V\in \Lambda$ a general member. 
We put $W:=f^{-1}(V)$. 
Then $D|_W$ is strongly horizontal 
with respect to $f_W:=f|_W:W\longrightarrow V$. 
\end{lem}

The following notion was introduced by M.~Reid. 

\begin{defn}[Nef and log big divisors]\label{log-big}
Let $(X,D)$ be lc and $L$ a $\bQ$-Cartier $\bQ$-divisor on $X$. 
The divisor $L$ is called {\em {nef and log big}} on $(X,D)$ 
if $L$ is nef and big, and $(L^{\dim W}\cdot W)>0$ 
for every center of log canonical 
singularities $W$ for the pair 
$(X,D)$. 
We note that an ample divisor is nef and log big. 
\end{defn}

We prepare some vanishing theorems. 
The following lemma is 
a special case of \cite[Lemma]{fk}, 
which is a variant of the Kawamata-Viehweg vanishing theorem 
(cf.~\cite{nori}). 

\begin{lem}\label{nori}
Let $X$ be a non-singular complete variety 
and $D$ a simple normal crossing divisor on $X$. 
Let $H$ be a nef and log big divisor on $X$. 
Then $H^i(X,K_X+D+H)=0$ for every $i>0$. 
\end{lem}

\begin{proof} 
If $D=0$, then $H^i(X,K_X+H)=0$ for every $i>0$ by 
the Kawamata-Viehweg vanishing theorem. 
So, we can assume that $D\ne 0$. 
Let $D_0$ be an irreducible component of $D$. 
We consider the following exact sequence: 
\begin{eqnarray*}
\cdots &\to & H^i(X,K_X+D-D_{0}+H)\to 
H^i(X,K_X+D+H)\\
&\to& H^i(D_0,K_{D_0}+{(D-D_0)}|_
{D_0}+H|_{D_0})\to \cdots.  
\end{eqnarray*}
By the inductions on the number of the irreducible components of 
$D$ and on $\dim X$, the first and the last terms are zero. 
Therefore, we obtain that $H^i(X,K_X+D+H)=0$ for every $i>0$. 
\end{proof}

The following proposition 
is a slight generalization of 
the Grauert-Riemenschneider vanishing theorem. 

\begin{prop}\label{relative}
Let $f:X\longrightarrow Y$ be a proper birational 
morphism from a non-singular variety $X$. 
Let $D$ be a simple normal crossing divisor on $X$, 
where $D$ may be zero. 
Assume that $f$ is an isomorphism at every generic point of 
center of log canonical singularities for the pair $(X,D)$. 
Then $R^if_*\omega _X(D)=0$ for $i>0$. 
\end{prop}
\begin{proof}
If $D=0$, then it is nothing but the Grauert-Riemenschneider vanishing 
theorem. So, we can assume that $D\ne 0$. 
Let $D_0$ be an irreducible component of $D$. 
We consider the following exact sequence: 
$$
\cdots \to R^if_*\omega_X(D-D_0)\to R^if_*\omega_X(D) \to 
R^if_*\omega_{D_0}((D-D_0)|_{D_0})\to \cdots. 
$$ 
Then $R^if_*\omega_X(D)=0$ for every positive $i$ by the 
same argument as in the proof of Lemma \ref{nori}. 
\end{proof}

The following corollary is an easy consequence of Proposition 
\ref{relative}, which will be used in Sections \ref{sec2} and 
\ref{sec3}. 

\begin{cor}\label{useful}
Let $f:Z\longrightarrow X$ be a proper birational 
morphism from a non-singular variety $Z$. 
Let $D$ be a reduced Weil divisor on $X$ such that 
$(X,D)$ is lc. 
Assume that $D'$ is the strict transform of $D$ and 
$f$ is an isomorphism over a Zariski open set $U$ of 
$X$ such that $U$ contains every generic point of 
center of log canonical singularities for 
the pair $(X,D)$. 
Then $f_*\omega_Z(D')\simeq \omega_X(D)$. 

We further assume that the $f$-exceptional locus 
$\Exc(f)$ $($see {\em{\ref{noanoa} (e)}} below$)$ 
and $D'\cup \Exc (f)$ are both simple normal 
crossing divisors on $Z$. 
Then $f$ is an isomorphism at every generic point of 
center of log canonical singularities for 
the pair $(Z,D')$ and $R^if_*\omega _Z(D')=0$ for $i>0$. 
In particular, $f$ induces a one to one correspondence 
between the generic points of center of log canonical 
singularities for the pair $(Z,D')$ and those for the pair 
$(X,D)$. 
\end{cor}
\begin{proof}
First, we write 
$$
K_Z+D'=f^*(K_X+D)+\sum _i a_iE_i, 
$$ 
where $E_i$ is an $f$-exceptional irreducible Cartier divisor 
on $Z$ for every $i$. 
Since $f$ is an isomorphism over $U$ 
that contains every generic point of log canonical singularities 
for the pair $(X,D)$, we have that $a_i>-1$ for 
every $i$. 
Thus, we obtain that $f_*\omega_Z(D')\simeq \omega_X(D)$. 

Next, let $W$ be a center of log canonical singularities for the 
pair $(Z,D')$. 
Then $W\not\subset \Exc(f)$ since 
$\Exc(f)$ and $D'\cup \Exc (f)$ are both simple normal 
crossing divisors. 
Therefore, $f$ is an isomorphism at every generic 
point of center of 
log canonical singularities for the pair $(Z,D')$. 
So, we can apply Proposition \ref{relative}. 
Thus, we obtain that $R^if_*\omega _Z(D')=0$ for $i>0$. 
The final statement is obvious by the above arguments. 
\end{proof}

\begin{rem}[Divisorial log terminal]\label{dlt}
The notion of {\em{dlt}} pairs may help the readers to 
understand Corollary \ref{useful}. 
Here, dlt is short for {\em{divisorial log terminal}}. 
It is one of the most useful variants of 
{\em{log terminal singularities}}. 
In this paper, however, we don't use it explicitly. 
So, we omit the details. 
For the precise definition and the basic properties 
of dlt pairs, see \cite{szabo}, \cite[\S 2.3]{kom}, 
and \cite{osamu}. 
\end{rem}

Finally, we fix the following notation and convention. 

\begin{say}[\bf{Notation and Convention}]\label{noanoa}
Let $\xN$ (resp.~$\mathbb Z_{\ge 0}$) be the set of positive
(resp.~non-negative) integers.
\begin{itemize}
\item[(a)] 
An {\em{algebraic fiber space}} $f:X\longrightarrow Y$ is a proper 
surjective morphism between non-singular projective 
varieties $X$ and $Y$ with connected fibers. 
\item[(b)] Let $f:X\longrightarrow Y$ be a dominant morphism 
between varieties. 
We put $\dim f:=\dim X-\dim Y$. 
\item[(c)] 
The words {\em{locally free sheaf}} and 
{\em{vector bundle}} are used interchangeably. 
\item[(d)] A Cartier (resp.~Weil) divisor $D$ 
on a normal variety $X$ and the associated line bundle 
(resp.~rank one reflexive sheaf) $\mathcal O_X(D)$ are 
used interchangeably if there is no danger of confusion. 
\item[(e)] Let $f:X\longrightarrow Y$ be a proper birational morphism 
between normal varieties. 
By the {\em{exceptional locus}} of $f$, we mean the subset 
$\{x\in X\ | \ \dim f^{-1}f(x)\geq 1\}$ of $X$, 
and denote it by $\Exc (f)$. 
We note that $\Exc(f)$ is of pure codimension one in 
$X$ if $f$ is birational and $Y$ is $\mathbb Q$-factorial. 
\item[(f)] When we use the desingularization theorem, 
we often forbid unnecessary blow-ups implicitly, 
that is, we don't use the weak Hironaka theorem but 
the original Hironaka theorem. 
Unnecessary blow-ups sometimes make the proof more 
difficult. 
We recommend the readers to see \cite[Remark 6-2 (iii)]{matsuki}, 
\cite[Resolution Lemma]{szabo}, and \cite[Corollary 7.9]{bev}, 
in particular, \cite[7.12 The motivation]{bev}. 
See also Remark \ref{resolution} below. 
My private note \cite{osamu} may help the readers to 
understand the subtleties of the desingularization theorem 
and various kinds of log terminal singularities.  
\item[(g)] Let $X$ be a normal variety and $D$ a $\mathbb Q$-divisor 
on $X$. 
A {\em{log resolution}} of $(X,D)$ is a proper birational morphism 
$g:Y\longrightarrow X$ such that $Y$ is non-singular, 
$\Exc (g)$ and $\Exc (g)\cup g^{-1}(\Supp D)$ are 
both simple normal crossing divisors. 
See \cite[Notation 0.4 (10)]{kom}. 
\end{itemize}
\end{say}

\section{Torsion-freeness and Vanishing theorem}\label{sec2}

In this section, we generalize Koll\'ar's torsion-free theorem 
and vanishing theorem: \cite[Theorem 2.1]{ko1}. 
The following is the main theorem of this section 
(see also Theorem \ref{hosoku} and Theorem \ref{vanishing}). 

\begin{theorem}\label{sawari}
Let $X$ be a non-singular projective variety, $D$ a 
simple normal crossing divisor on $X$, $Y$ 
an arbitrary $($reduced$)$ projective variety, and 
$f:X\longrightarrow Y$ a surjective morphism. 
Then 
\begin{itemize}
\item[(i)] If $H$ is an ample divisor on $Y$, then 
$H^j(Y,H\otimes R^if_*\omega_{X/Y}(D))=0$ for $j>0$. 
\end{itemize}
Assume furthermore 
that $D$ is strongly horizontal with respect to $f$. 
Then 
\begin{itemize}
\item[(ii)] $R^if_*\omega_{X}(D)$ is torsion-free for $i\geq 0$. 
\item[(iii)]$R^if_*\omega_{X}(D)=0$ if $i>\dim f$. 
Furthermore, if $D\ne 0$, then $R^if_*\omega_{X}(D)=0$ 
for $i\geq \dim f$.    
\end{itemize}
\end{theorem}

The statement (iii) is obvious by (ii). 
So, it is sufficient to prove (i) and (ii). 
In Theorem \ref{torsionfree}, we treat (ii) in a more general 
setting. We will prove (i), which is not 
used later, in the subsection \ref{subsec2.2}.

\subsection{Torsion-free theorem}\label{subsec2.1}
The following is the main theorem of this subsection. 
It is a special case of \cite[Theorem 3.2 (i)]{am}. 
We adopt Arapura's proof of torsion-freeness (see the proof 
of Theorem 1 in \cite{a}). 
This proof is suitable for our paper.  

\begin{thm}[Torsion-freeness]\label{torsionfree}
Let $f:X\longrightarrow Y$ be a surjective morphism from a non-singular 
projective variety $X$ 
to a {\em{(}}possibly singular{\em{)}} projective variety $Y$. 
Let $D$ be a simple normal crossing divisor on $X$. 
Assume that $D$ is strongly horizontal with respect to $f$. 
Let $L$ be a semi-ample line bundle on $X$. 
Then, for all $i$, 
the sheaves $R^if_*(\omega_X(D)\otimes L)$ are 
torsion-free. 
In particular, $R^if_*\omega _X(D)$ is torsion-free for 
every $i$. 
\end{thm}
\begin{proof}
In order to explain the plan of the proof, let us 
introduce the following notation, where $f:X\longrightarrow Y$ 
and the divisor $D$ are as in the statement. 

\begin{tabular}{rcl}
$P_n^{\log}$ & : & $\begin{cases}
\begin{minipage}{9.5cm} 
If $\dim \Supp (R^if_*\omega_X(D)_{\tor})\leq 
n$ for all $i$, 
then the sheaves $R^if_*\omega_{X}(D)$ are torsion-free for 
all $i$, 
where $R^if_*\omega_X(D)_{\tor}$ 
is the torsion part of $R^if_*\omega_X(D)$.
\end{minipage}\end{cases}$\\
$P^{\log}(Y)$ & : & $\begin{cases}\begin{minipage}{9.5cm} 
The sheaf $R^if_*\omega_X(D)$ is 
torsion-free for every $i$. \end{minipage}\end{cases}$\\
$Q^{\log}(Y)$ & : & $\begin{cases}\begin{minipage}{9.5cm}
Then the sheaf $R^if_*(\omega_X(D)\otimes L)$ is 
torsion-free for every $i$. \end{minipage}\end{cases}$
\end{tabular}
It is enough to prove the following four claims: 
\begin{itemize}
\item $P^{\log}(Y)$ implies $Q^{\log}(Y)$.  

\item $P^{\log}(Y)$ when $Y$ is a curve. 

\item $Q^{\log}(\mathbb P^1)$ implies $P_0^{\log}$. 

\item $P_{n-1}^{\log}$ implies $P_n^{\log}$.   
\end{itemize}
\begin{step}[$P^{\log}(Y)$ implies $Q^{\log}(Y)$]\label{step1} 
Since $L$ is semi-ample, there exists an $m>0$ for which 
$L^m$ is generated by global sections. 
Hence, by Bertini's theorem, we can find $B\in |L^m|$ such that 
$B$ is smooth and $B+D$ is a simple normal crossing divisor on $X$. 
Let $\pi:Z\longrightarrow X$ be the $m$-fold cyclic covering 
branched along $B$. 
Then $\pi_*\omega_Z\simeq \bigoplus _{j=0}^{m-1}
\omega_X\otimes L^j$. 
Therefore, $\omega_X(D)\otimes L$ is a direct summand of $\pi_*
\omega_{Z}(\pi^*D)$. 
Since $\pi$ is finite, 
we have 
$$
R^i(f\circ \pi)_*\omega_Z({\pi^*D})=R^if_*(\pi_*\omega_Z(\pi^* D))
=\bigoplus_{j=0}^{m-1}R^if_*(\omega_X(D)\otimes L^j). 
$$
By $P^{\log}(Y)$, the left hand side is torsion-free. 
We note that $\pi^*D$ is a simple normal crossing divisor 
on $Z$ and strongly horizontal with 
respect to $f\circ \pi$. 
Then, so is $R^if_*(\omega_X(D)\otimes L)$. 
\end{step}
\begin{step}[$P^{\log}(Y)$ when $Y$ is a curve] 
Now suppose that $Y$ is a curve. 
By the Stein factorization, we can assume that 
$Y$ is smooth. 
Let $Y_0$ be a non-empty Zariski open set of $Y$ such 
that $f$ is smooth and $D$ is relatively normal crossing 
over $Y_0$. 
By blowing up $X$, we can assume that 
$\Supp (f^{-1}P\cup D)$ is simple normal crossing for 
$P\in Y\setminus Y_0$ (cf.~Corollary \ref{useful}). 
If $f:X\longrightarrow Y$ is semi-stable, then 
the theorem follows from \cite[(5.7)]{sz} (see also 
Theorem \ref{vmhs} and 
Step \ref{step4} in the proof of Theorem \ref{can-ext} 
below). 
If $f:X\longrightarrow Y$ is not semi-stable, 
then we apply the semi-stable reduction theorem (cf.
~\cite[Chapter II]{kmms} and 
~\cite[I.9]{ssu2}). 
We consider the following commutative diagram: 
$$
\begin{CD}
X@<{\alpha}<<X' @<{\beta}<<{\widetilde X}\\ 
@V{f}VV @V{f'}VV @VV{\widetilde f}V \\
Y@<<{\tau}<Y'@=Y',  
\end{CD}
$$ 
where $\tau:Y'\to Y$ is a finite morphism from a non-singular 
projective curve and $X'$ is the normalization of 
$X\times _{Y}Y'$, and $\beta$ is a birational morphism 
such that $\widetilde f:\widetilde X\to Y'$ is semi-stable. 
We note that we can assume that $\beta$ is an isomorphism over a 
non-empty Zariski open set of $Y'$ (cf.~\ref{noanoa} (f)). 
Then $R^if_*\omega_X(D)$ is a direct summand of 
$\tau_*R^i{f'}_*\omega_{X'}(\alpha ^*D)$. 
So, it is sufficient to check the local freeness of 
$R^i{f'}_*\omega_{X'}(\alpha ^*D)$. 
We note that $R^i{f'}_*\omega_{{X'}/Y'}(\alpha ^*D)\simeq 
R^i{\widetilde f}_*\omega_{{\widetilde X}/Y'}(D')$, 
where $D'$ is the proper transform of $\alpha ^*D$ 
(cf.~Corollary \ref{useful}). 
We can assume that $\Supp ({\widetilde f}^{-1}P\cup D')$ 
is simple normal crossing for 
every $P\in Y'$ (cf.~\cite[I.9]{ssu2}). 
Since $\widetilde f$ is semi-stable, 
$R^i{\widetilde f}_*\omega_{{\widetilde X}/Y'}(D')$ is 
locally free by the above argument. 
Thus, we obtain that $R^if_*\omega_X(D)$ is locally free 
for every $i$. 
\end{step}
\begin{step}[$Q^{\log}(\mathbb P^1)$ implies $P_0^{\log}$] 
We assume that 
the sheaf $R^if_*\omega_X(D)$ has 
torsion supported on a finite set of 
points $S:=\{p_1, \cdots , p_r\}$ for some $i$. 
Now, take a pencil of hyperplane sections of $Y$. 
We can assume that the base locus is disjoint from $S$ and 
that the preimage of the base locus in $X$ is smooth 
and meets all the centers of log canonical singularities 
of $(X,D)$ transversally. 
Blow up the base locus and its preimage in $X$ to get a diagram. 
$$
\begin{matrix}
f':&X' \ \ \ & \longrightarrow  & \ \ \ Y' \\
&{\searrow}&  &  {\swarrow g} \\
& & \mathbb P^1 &
\end{matrix}
$$
Let $\mathcal H$ be an ample line bundle on $Y'$. 
Replacing $\mathcal H$ by $\mathcal H^k$ for some $k\gg 0$, 
we can assume that 
$R^pg_{*}(\mathcal H\otimes R^qf'_*\omega_{X'}(D'))=0$ 
for all $p>0$ and for all $q$, 
where $D'$ is the strict transform of $D$ on $X'$. 
Therefore, the spectral sequence collapses to give isomorphisms 
$g_*(\mathcal H\otimes R^qf'_*\omega_{X'}(D'))=R^q(g\circ f')_*
({f'}^*\mathcal H\otimes \omega_{X'}(D'))$. 
By $Q^{\log}(\mathbb P^1)$, 
the right hand side is torsion-free. However, by the 
assumption, the sheaf $\mathcal H\otimes R^if'_*\omega_{X'}(D')$ has 
torsion supported at the points $p_j\in Y'$. 
Therefore, $g_*(\mathcal H\otimes R^if'_*\omega_{X'}(D'))$ 
has torsion at the 
points $g(p_j)$. 
This is a contradiction. 
Thus, the sheaf $R^if_*\omega_X(D)$ must be torsion-free. 
\end{step}
\begin{step}[$P_{n-1}^{\log}$ implies $P_n^{\log}$]   
Assume that $\dim \Supp (R^if_*\omega_X(D)_{\tor})\leq 
n$ for all $i$. 
We suppose that for some $i$ the sheaf $R^if_*\omega_X(D)$ is not 
torsion-free. 
Then there must be a positive dimensional component of 
$\Supp (R^if_*\omega_X(D)_{\tor})$ by $P^{\log}_{n-1}$. 

Let $\mathcal H$ be a very ample line bundle on $Y$ and 
let $B\in |\mathcal H|$ be a general member such that 
$f^*B$ is smooth and $f^*B+D$ is a simple normal crossing divisor. 
Then $R^if_*\omega_{f^*B}(D|_{f^*B})\simeq 
R^if_*\omega_X(D)\otimes \mathcal O_B(B)$. 
Applying $P^{\log}_{n-1}$ to $(f^*B, D|_{f^*B})\longrightarrow B$, 
we obtain that the left hand side is torsion-free. 
This contradicts the assumption that $R^if_*\omega_X(D)$ has 
torsion. 
So, we obtain the required result. 
\end{step}
Therefore, we complete the proof. 
\end{proof}

We can omit the assumption that $X$ and $Y$ are projective in 
Theorem \ref{torsionfree}. 
We will use Theorem \ref{hosoku} in the proof of 
Theorem \ref{can-ext}. 

\begin{thm}\label{hosoku}
Let $f:X\longrightarrow Y$ be a projective surjective 
morphism from a non-singular variety 
to a {\em{(}}possibly singular{\em{)}} variety. 
Let $D$ be a simple normal crossing divisor on $X$. 
Assume that $D$ is strongly horizontal with respect to $f$. 
Let $L$ be a semi-ample line bundle on $X$. 
Then, for all $i$, 
the sheaves $R^if_*(\omega_X(D)\otimes L)$ are 
torsion-free. 
In particular, $R^if_*\omega _X(D)$ is torsion-free for 
every $i$. 
\end{thm}
\begin{proof}
By Step \ref{step1} in the proof of Theorem \ref{torsionfree}, 
it is enough to prove that 
$R^if_*\omega _X(D)$ is torsion-free for 
every $i$. 
Since the statement is local, we can 
shrink $Y$ and assume that $Y$ is quasi-projective. 
We can take a suitable compactification and assume that 
$X$ and $Y$ are both projective (see Remark 
\ref{resolution} below). 
Thus, by Theorem \ref{torsionfree}, we obtain the 
required result. 
\end{proof}

\begin{rem}\label{resolution}
Here, we had better use 
Szab\'o's resolution lemma: \cite[Resolution Lemma]{szabo}. 
See also \cite[Remark 6-2 (iii)]{matsuki} 
or \cite[Corollary 7.9]{bev}. 
We recommend the readers to see 
\cite[7.12 The motivation]{bev}. 
\end{rem}

The following example implies that 
if $D$ is not strongly horizontal, 
then the torsion-freeness is not necessarily 
true. 

\begin{ex}\label{takagi}
Let $Y$ be a non-singular projective surface 
and $f:X\longrightarrow Y$ be a blow-up at a point $P\in Y$. 
We put $D:=f^{-1}(P)$. 
We consider the following exact sequence: 
$$
0\longrightarrow \omega_X\longrightarrow 
\omega_X(D)\longrightarrow \omega_D\longrightarrow 0. 
$$ 
Then we obtain that $R^1f_*\omega_X(D)\simeq R^1f_*\omega_D\simeq 
H^1(D,\omega_D)\simeq \mathbb C_P$ by 
the Grauert-Riemenschneider vanishing theorem. 
So, $R^1f_*\omega_X(D)$ is not torsion-free. 
\end{ex}

\begin{cor}\label{norimatsu}
Let $f:X\longrightarrow Y$ and $D$ be as in {\em{Theorem \ref{hosoku}}}. 
Assume that $L$ is a relatively ample line bundle on $X$. 
Then $R^if_*(\omega_X(D)\otimes L)=0$ for $i>0$. 
\end{cor}
\begin{proof}
It is essentially the same as \cite[Corollary 2 (i)]{a}. 
It is sufficient to use Lemma \ref{nori}, 
instead of the Kodaira vanishing theorem.  
\end{proof}

\subsection{Vanishing theorem}\label{subsec2.2}

The following is a slight generalization of Koll\'ar's vanishing 
theorem: \cite[Theorem 2.1 (iii)]{ko1}. 
It is also a special case of \cite[Theorem 3.2 (ii)]{am}. 
We adopt 
the proof of \cite[9.14 Theorem]{ks}. 
We will not use this vanishing theorem later. 
So, the readers can skip this subsection. 

\begin{thm}[Vanishing Theorem]\label{vanishing}
Let $f:X\longrightarrow Y$ be a morphism from a non-singular 
projective variety 
$X$ onto a variety $Y$. 
Let $D$ be a simple normal crossing divisor on $X$. 
Let $H$ be an ample line bundle on $Y$. 
Then 
$$
H^i(Y,H\otimes R^jf_*\omega _X(D))=0 
$$ 
for $i>0$ and  $j\geq 0$.   
\end{thm}

\begin{proof}
Let $n$ be a positive integer such that $n\geq 2$ and 
the linear system $|H^n|$ is base point free. 
Take a general member $E\in |H^n|$ such that 
$Z:=f^{-1}(E)$ is smooth and $Z\cup D$ is a simple normal crossing 
divisor. 
By \cite[5.1 b)]{ev}, 
\begin{equation}\label{injective}
H^i(X,\omega_X(D)\otimes f^{*}H)\longrightarrow 
H^i(X,\omega_X(D)\otimes f^{*}H^{1+kn})
\end{equation} 
is injective for $i\geq 0$. 
We prove the theorem by induction on $\dim Y$. 
The assertion is evident if $\dim Y=0$. 
We have an exact sequence: 
$$
0\to \omega_X(D)\otimes f^*H^t\to 
\omega_X(D)\otimes f^*H^{t+n}\to 
\omega_Z(D|_Z)\otimes (f^*H^t|_Z)\to 0.  
$$ 
Using induction and the corresponding long exact sequence, 
we obtain that 
$$
H^i(Y, H^t\otimes R^jf_*\omega_X(D))\simeq 
H^i(Y, H^{t+n}\otimes R^jf_*\omega_X(D))  
$$ 
for $i\geq 2$. 
By the Serre vanishing theorem, 
$$
H^i(Y, H^{t+kn}\otimes R^jf_*\omega_X(D))=0 
$$ 
for $k\gg 0$. 
Thus, 
$$
H^i(Y, H^{t}\otimes R^jf_*\omega_X(D))=0 
$$ 
for $t\geq 1$ and $i\geq 2$. 
Once this much of the theorem is established, the Leray 
spectral sequence 
$$
E^{p,q}_2=H^p(Y, H^t\otimes R^q f_*\omega_X(D))\Longrightarrow 
E^{p+q}=H^{p+q}(X, \omega_X(D)\otimes f^*H^t) 
$$ 
has only two columns, 
and therefore it degenerates. 
This means that 
$$
0\longrightarrow E^{1,j}_2\longrightarrow E^{j+1}\longrightarrow 
E^{0,j+1}_2\longrightarrow 0. 
$$ 
Then we have 
$$
\begin{CD}
0 &\longrightarrow & H^1(Y, H\otimes R^if_*\omega_X(D))&
\longrightarrow &
H^{j+1}(X, \omega_X(D)\otimes f^*H)\\ 
& &\downarrow & & \downarrow\\ 
0 &\longrightarrow & H^1(Y, H^{1+kn}\otimes R^if_*\omega_X(D))&
\longrightarrow &
H^{j+1}(X, \omega_X(D)\otimes f^*H^{1+kn}). 
\end{CD}
$$ 
Using (\ref{injective}), this implies that 
$$
H^1(Y, H\otimes R^jf_*\omega_X(D))\longrightarrow 
H^1(Y, H^{1+kn}\otimes R^jf_*\omega_X(D)) 
$$ 
is injective for every $k$. 
As before, by the Serre vanishing theorem, 
this implies that $H^1(Y, H\otimes R^jf_*\omega_X(D))=0$. 
We complete the proof. 
\end{proof}

\section{Variation of mixed Hodge structures}\label{sec3}

To investigate $R^if_*\omega_X(D)$, we use the notion of 
gradedly polarized variation of mixed Hodge structures. 
We note that 
all the variations of mixed Hodge structures which we treat in 
this section are geometric.  

\subsection{Canonical Extension}\label{subsec3.1}

In this subsection, we generalize \cite[Theorem 2.6]{ko2} and 
\cite[Theorem 1]{n}. 

\begin{say}\label{setting} 
Let $f:X\longrightarrow Y$ be a projective surjective morphism 
between non-singular varieties over $\mathbb C$. 
Let $D$ be a simple normal crossing divisor on $X$ such 
that $D$ is strongly horizontal. 
Assume that there is a non-empty Zariski open set $Y_0$ 
of $Y$ such that $\Sigma:=Y\setminus Y_0$ is a simple 
normal crossing divisor on $Y$ 
and that $f_0:X_0\longrightarrow Y_0$ 
is smooth and $D_0$ is 
relatively normal crossing over 
$Y_0$, where $X_0:=f^{-1}(Y_0)$, 
$f_0:=f|_{X_0}$ and $D_0:=D\cap X_0$. 
The local system $R^if_{0*}\mathbb C_{X_0- D_0}$ on 
$Y_0$ forms a gradedly polarized 
variation of mixed Hodge structure (see \cite{ssu}). 
\end{say}

\begin{say}\label{setting2} 
Assume that all the local monodromies of the 
local system $R^kf_{0*}\mathbb C_{X_0-D_0}$ 
around $\Sigma$ are unipotent. 
Put $\mathcal H^k_0:=(R^kf_{0*}\mathbb C_{X_0-D_0})\otimes 
\mathcal O_{Y_0}$ and let $F^p(\mathcal H^k_0)$ 
be the $p$-th Hodge filtration of $\mathcal H^k_0$. 
Let $\mathcal H^k_Y$ be 
the canonical extension of $\mathcal H^k_0$ 
to $Y$. 
Then there exists an extension 
$F^p(\mathcal H^k_Y)$ of $F^p(\mathcal H^k_0)$, 
which is locally free. See \cite[\S 5 The geometric case]{sz},
\cite[I.10]{ssu2}, and 
\cite[Lemma 1.11.2]{kashi}. 
We note that $F^p(\mathcal H^k_Y)=j_*F^p(\mathcal H^k_0) \cap 
\mathcal H_Y^k$, where $j:Y_0\longrightarrow Y$ is 
the natural inclusion. 
As stated above, in 
this paper, we only treat {\em{geometric}} 
gradedly polarized variation of 
mixed Hodge structures. 
\end{say}

The following is the main theorem of this subsection (see also 
Theorem \ref{loc-free}). 
The proof is essentially the same as the proof of \cite
[Theorem 1]{n}. 

\begin{thm}\label{can-ext}
Under the same notation as in {\em{\ref{setting}}}, 
let $\omega_{X/Y}:=\omega_X\otimes f^*\omega^{-1}_Y$ and 
$d:=\dim f$. 
Assume that all the local monodromies of the 
local system $R^{d+i}f_{0*}\mathbb C_{X_0-D_0}$ 
around $\Sigma$ are unipotent. 
Then there exists an isomorphism 
$$
R^if_*\omega_{X/Y}(D)\simeq F^d(\mathcal H^{d+i}_Y).  
$$ 
In particular, 
$R^if_*\omega_{X/Y}(D)$ is locally free. 
\end{thm}

Our proof of Theorem \ref{can-ext} relies on the following theorem. 
We can take it out from (\cite[\S 5 The geometric case]{sz}) 
with a little effort. 
See also \cite[I.10]{ssu2}. 

\begin{thm}[{\cite[\S 5]{sz}}]\label{vmhs}
Let $f:X\longrightarrow Y$ be a projective surjective morphism 
from a non-singular variety $X$ onto a non-singular 
curve $Y$. 
Let $D$ be a simple normal crossing divisor on $X$. 
Assume that there is a divisor $\Sigma$ on $Y$ such that 
$f$ is smooth over $Y_0:=Y\setminus \Sigma$ and 
$D$ is relatively normal crossing over $Y_0$ and 
that $C\cup D$ is a simple normal crossing divisor, 
where $C:=(f^*\Sigma)_{red}$. 
Assume that all the local monodromies on $R^if_{0*}\mathbb 
C_{X_0-D_0}$ around $\Sigma$ are unipotent. 
Then 
$$
\mathcal H^i_Y \simeq \mathbf R^if_*\Omega_{X/Y}^{\bullet} 
(\xlog(C\cup D)) 
$$ 
and 
$$
F^p(\mathcal H^i_Y)\simeq \mathbf R^if_*F^p(\Omega_{X/Y}^{\bullet} 
(\xlog(C\cup D))) 
$$ 
for all $p$. 

Here, $\Omega^{\bullet}_{X/Y}(\log (C\cup D))$ is 
the {\em{relative log complex}}{\em{:}} 
$$
\Omega^{\bullet}_{X/Y}(\log (C\cup D)):=
\Omega^{\bullet}_{X}(\log (C\cup D))/f^*\Omega^1_Y(\log \Sigma)\wedge 
\Omega^{\bullet}_{X}(\log (C\cup D))[-1], 
$$ 
and $K^{\bullet}=F^p(\Omega^{\bullet}_{X/Y}(\log (C\cup D)))$ 
is a complex such that 
$$
K^q=
\begin{cases}
0 \hspace{40mm}\text{if} \quad q<p\\
\Omega^{q}_{X/Y}(\log (C\cup D)) \hspace{10mm}\text{otherwise}. 
\end{cases}
$$ 
\end{thm}
\begin{proof}[{Proof of {\em{Theorem \ref{can-ext}}}}] 
By Corollary \ref{useful} 
and \ref{noanoa} (f), we can assume that 
$C\cup D$ is a simple normal crossing divisor on 
$X$ without loss of generality, where $C:=(f^*\Sigma)_{red}$. 
\setcounter{step}{0}
\begin{step}[The case when $\dim Y=1$]\label{step4}
By Theorem \ref{vmhs}, we have 
$$
F^d(\mathcal H^{d+i}_Y)\simeq \mathbf 
R^if_*\Omega_{X/Y}^{d} 
(\xlog(C\cup D)).  
$$ 
On the other hand, $\Omega_{X/Y}^{d} 
(\xlog(C\cup D))\simeq \omega_{X/Y}(C-f^*\Sigma+D)$. 
Therefore, if $f$ is semi-stable, 
then $R^if_*\omega_{X/Y}(D)\simeq F^d(\mathcal H^{d+i}_Y)$. 
If $f$ is not semi-stable, then we use 
the covering 
arguments in \cite[Lemma 2.11]{ko2} and \cite[Lemma 1.9.1]{kashi}. 
Thus, we obtain 
$F^d(\mathcal H^{d+i}_Y)
\simeq 
R^if_*(\Omega_{X/Y}^{d} 
(\xlog(C\cup D))\otimes \mathcal O_X(\sum (a_i-1)C_i)) 
\simeq R^if_*\omega_{X/Y}(D)$, where $f^*\Sigma:=\sum a_i C_i$. 
Note that the middle term is 
the {\em{upper canonical extension}} (cf.~\cite[Definition 2.3]{ko2}) 
and 
there is no difference between the 
{\em{canonical extension}} 
and the upper canonical extension (the {\em{right canonical extension}} 
in the proof of \cite[Lemma 1.9.1]{kashi}), 
since all the local monodromies 
are unipotent. See Remark \ref{dasoku} below. 
\end{step}

\begin{step}
[The case when $l:=\dim Y\geq 2$]\label{iran} 
We shall prove by induction on $l$. 

By Step\ref{step4}, there is an open subset $Y_1$ of $Y$ 
such that $\codim (Y\setminus Y_1)\geq 2$ and 
that 
$$
R^if_*\omega_{X/Y}(D)|_{Y_1}\simeq F^d(\mathcal H_Y^{d+i})|_{Y_1}. 
$$ 
Since $F^d(\mathcal H_Y)$ is locally free, we obtain a 
homomorphism 
$$
\varphi_Y^i: 
R^if_*\omega_{X/Y}(D)\longrightarrow F^d(\mathcal H_Y^{d+i}). 
$$ 
By Theorem \ref{hosoku}, $\Ker \varphi_Y^i=0$. 
We put $G_Y^i:= \Coker \varphi_Y^i$. 
Taking a general hyperplane cut, we see that $\Supp G_Y^i$ 
is a finite set by the induction hypothesis. 
Assume that $G_Y^i\ne 0$. 
Take a point $P\in G_Y^i$. 
Let $\mu:W\longrightarrow Y$ be the blowing up at $P$ and 
set $E=\mu^{-1}(P)$. 
Then $E\simeq \mathbb P^{l-1}$. 
By \ref{noanoa} (f), we can take a projective birational 
morphism $p:X'\longrightarrow X$ from a non-singular 
variety $X'$ with the following properties: 
\begin{itemize}
\item[(i)] the composition 
$X'\longrightarrow X\longrightarrow Y\dashrightarrow W$ 
is a morphism. 
\item[(ii)] $p$ is an isomorphism over $X_0$. 
\item[(iii)] $\Exc (p)$ and $\Exc (p)\cup D'$ are 
simple normal crossing divisors on $X'$, 
where $D'$ is the strict transform of $D$. 
\end{itemize}
By Corollary \ref{useful}, 
we obtain that $R^if_*\omega _{X/Y}(D)\simeq 
R^i(f\circ p)_*\omega_{X'/Y}(D')$ for every $i$. 
We note that $D'$ is strongly horizontal with 
respect to $f\circ p$.
 By replacing $(X,D)$ with 
$(X',D')$, we can assume that there is a morphism 
$g:X\longrightarrow W$ such that $f=\mu\circ g$. 
Since $g:X\longrightarrow W$ is in the same situation as $f$, 
we obtain the exact sequence: 
$$
0\to R^ig_*\omega_{X/W}(D)\to F^d(\mathcal H_W^{d+i})\to 
G_W^i \to 0.
$$ 
Tensoring $\mathcal O_W(\nu E)$ for $0\leq \nu \leq l-1$ and 
applying $R^j\mu_*$ for $j\geq 0$ to each $\nu$, 
we have a exact sequence 
\begin{eqnarray*}
0 & \to & \mu_*(R^ig_*\omega_{X/W}(D)\otimes \mathcal O_W(\nu E)) 
\to \mu_*(F^d(\mathcal H_W^{d+i})\otimes \mathcal O_W(\nu E)) \\ 
&\to &\mu_*(G_W^i\otimes \mathcal O_W(\nu E)) 
\to  R^1\mu_*(R^ig_*\omega _{X/W}(D)\otimes \mathcal O_{W}(\nu E)) \\ 
&\to & R^1\mu_*(F^d(\mathcal H_W^{d+i})\otimes \mathcal O_W(\nu E)) 
\to 0 
\end{eqnarray*} 
and $R^q\mu_*(R^ig_*\omega_{X/W}(D)\otimes \mathcal O_W(\nu E))
\simeq R^q\mu_*(F^d(\mathcal H_W^{d+i})\otimes \mathcal O_W(\nu E))
$ for $q\geq 2$. 

By Lemma \ref{pull} below, $F^d(\mathcal H_W^{d+i})\simeq 
\mu^*F^d(\mathcal H_Y^{d+i})$. 
We have 
$$
\mu_*(F^d(\mathcal H_W^{d+i})\otimes 
\mathcal O_W(\nu E))\simeq F^d(\mathcal H_Y^{d+i})
$$ 
and 
$$
R^q\mu_*(F^d(\mathcal H_W^{d+i})\otimes 
\mathcal O_W(\nu E))=0
$$ for $q\geq 1$. 
Therefore, $R^q\mu_*(R^ig_*\omega_{X/W}(D)\otimes \mathcal O_W(\nu E))=0$ 
for $q\geq 2$ and 
\begin{eqnarray*}
0 & \to & \mu_*(R^ig_*\omega_{X/W}(D)\otimes \mathcal O_W(\nu E)) 
\to \mu_*(F^d(\mathcal H_W^{d+i})\otimes \mathcal O_W(\nu E)) \\ 
&\to &\mu_*(G_W^i\otimes \mathcal O_W(\nu E)) 
\to  R^1\mu_*(R^ig_*\omega _{X/W}(D)\otimes \mathcal O_{W}(\nu E))\to 0 
\end{eqnarray*} 
is exact. 
Since $\omega_W=\mu^*\omega_Y\otimes \mathcal O_W((l-1)E)$, 
we have a spectral sequence 
$$
E^{p,q}_2=R^p\mu_*(R^qg_*\omega_{X/W}(D)\otimes \mathcal O_W((l-1)E)) 
\Longrightarrow R^{p+q}f_*\omega_{X/Y}(D). 
$$ 
However, $E^{p,q}_2=0$ for $p\geq 2$ by the above argument; 
thus 
\begin{eqnarray*}
0&\to& R^1\mu_*R^{i-1}g_*\omega_{X/Y}(D)\to R^if_*\omega_{X/Y}(D)\\ 
&\to& 
\mu_*(R^ig_*\omega _{X/W}(D)\otimes \mathcal O_W((l-1)E))\to 0. 
\end{eqnarray*} 
By Theorem \ref{hosoku}, 
$R^1\mu_*R^{i-1}g_*\omega_{X/Y}(D)=0$. 
Therefore, for $q\geq 1$, we 
obtain 
\begin{itemize}
\item[(a)] $R^if_*\omega_{X/Y}(D)\simeq \mu_*(R^ig_*\omega_{X/W}(D)
\otimes \mathcal O_W((l-1)E))$ and  
\item[(b)] $R^q\mu_*(R^ig_*\omega_{X/W}(D)\otimes \mathcal O_W((l-1)E))=0$. 
\end{itemize} 
Next, we shall consider the following commutative 
diagram: 

$$
\begin{matrix}
0 & &0\\
\downarrow & &\downarrow\\ 
R^ig_*\omega _{X/W}(D)\otimes \mathcal O_W((l-2)E)&
\longrightarrow &
R^ig_*\omega _{X/W}(D)\otimes \mathcal O_W((l-1)E) \\ 
\downarrow & & \downarrow\\ 
F^d(\mathcal H_W^{d+i})\otimes \mathcal O_W((l-2)E)
& \longrightarrow &
F^d(\mathcal H_W^{d+i})\otimes \mathcal O_W((l-1)E) \\ 
\downarrow & & \downarrow\\
G_W^{i}\otimes \mathcal O_W((l-2)E)
&\longrightarrow & 
G_W^{i}\otimes \mathcal O_W((l-1)E)\\ 
\downarrow & & \downarrow\\ 
0 &  &0
\end{matrix}
$$ 
By applying $\mu_*$, we have 
$$
\begin{matrix}
0 & &0\\
\downarrow & &\downarrow\\ 
\mu_*(R^ig_*\omega _{X/W}(D)\otimes \mathcal O_W((l-2)E))
&\longrightarrow &
\mu_*(R^ig_*\omega _{X/W}(D)\otimes \mathcal O_W((l-1)E)) \\
\downarrow & & \downarrow\\ 
F^d(\mathcal H_Y^{d+i})
& \simeq &
F^d(\mathcal H_Y^{d+i})\\
\downarrow & & \downarrow\\
\mu_*(G_W^{i}\otimes \mathcal O_W((l-2)E))
&\longrightarrow & 
\mu_*(G_W^{i}\otimes \mathcal O_W((l-1)E))\\
& & \downarrow\\ 
&  &0
\end{matrix}
$$ 
By (a) and (b), 
$G_Y^i\simeq \mu_*(G_W^i\otimes \mathcal O_W((l-1)E))$ 
and 
$$
\mu_*(G_W^i\otimes \mathcal O_W((l-2)E))\to 
\mu_*(G_W^i\otimes \mathcal O_W((l-1)E)) 
$$ 
is surjective. 
Since $\dim \Supp G_W^i=0$ and $E\cap \Supp G_W^i\ne \emptyset$, 
it follows that $G^i_W=0$ by Nakayama's lemma. 
Therefore, $G_Y^i=0$. 
This prove the theorem. 
\end{step}
\end{proof}

The following lemma played an essential role in the proof 
of Theorem \ref{can-ext}. 

\begin{lem}\label{pull}
Let $f:X\longrightarrow Y$ be and $D$ be as in {\em{\ref{setting}}}. 
Let $\pi:V\longrightarrow Y$ be 
a morphism from a non-singular variety such 
that $\Supp \pi^{-1}(\Sigma)$ is a simple normal crossing divisor on $V$. 
Then we obtain the gradedly polarized variation of 
mixed Hodge structures on $V_0:=V\setminus \pi^{-1}(\Sigma)$ 
by the base change. 
Assume that all the local monodromies of the local system 
$R^kf_{0*}\mathbb C_{X_0-D_0}$ around $\Sigma$ are unipotent. 
Then $\pi^*F^p(\mathcal H_Y^k)\simeq F^p(\mathcal H_V^k)$, 
where $\mathcal H_V^k$ is the canonical extension of $\mathcal 
H_{V_0}:=\pi^*
\mathcal H_0^k$ to $V$. 
\end{lem}
\begin{proof}
Note that $F^p(\mathcal H_Y^k)=j_*F^p(\mathcal H_0^k)\cap 
\mathcal H_Y^k$, where $j:Y_0\longrightarrow Y$ 
is the natural inclusion. See, 
for example, \cite[I.10]{ssu2} and \cite[Lemma 1.11.2]{kashi}. 
Then we have $\pi^*F^p(\mathcal H_Y^k)\simeq F^p(\mathcal H_V^k)$. 
See \cite[Proposition 1]{ka2}. 
\end{proof}

By using the unipotent reduction theorem, we obtain the 
following theorem. 

\begin{thm}\label{loc-free}
We use the same notation and assumptions as in {\em{\ref{setting}}}. 
We put $\omega _{X/Y}:=\omega_X\otimes \omega_Y^{-1}$ 
and $d:=\dim f$. 
Then $R^if_*\omega_{X/Y}(D)$ is locally free. 
More precisely, $R^if_*\omega_{X/Y}(D)$ is 
the upper canonical extension of $R^{i}{f_0}_*\omega_{X_0/Y_0}(D_0)$ 
{\em{(}}see \cite[Definition 2.3]{ko2} and Remark 
{\em{\ref{dasoku}}} below{\em{)}}. 
\end{thm}
\begin{proof}(cf.~{\cite[Reduction 2.12]{ko2}}) 
It is sufficient to prove the local freeness of 
$R^if_*\omega_{X/Y}(D)$. 
We already checked that $R^if_*\omega_{X/Y}(D)$ is the upper 
canonical extension in codimension one (see Step \ref{step4} 
in the proof of Theorem \ref{can-ext}). 
We use the unipotent reduction theorem with 
respect to the local system $R^{d+i}f_{0*}\mathbb C_{X_0-D_0}$. 
First, we can assume that $\Supp (D\cup f^{-1}(\Sigma))$ is a 
simple normal crossing divisor (cf.~Corollary \ref{useful} and 
\ref{noanoa} (f)). 
We consider the following commutative diagram: 
$$
\begin{CD}
X@<{\alpha}<<X' @<{\beta}<<{\widetilde X}\\ 
@V{f}VV @V{f'}VV @VV{\widetilde f}V \\
Y@<<{\tau}<Y'@=Y',  
\end{CD}
$$ 
where $\tau:Y'\to Y$ is a finite morphism from a non-singular 
variety obtained by Kawamata's covering trick, 
$X'$ is the normalization of 
$X\times _{Y}Y'$, $\beta$ is a projective birational morphism 
from a non-singular variety $\widetilde X$, 
and $D'$ is a simple normal crossing divisor 
on $\widetilde X$ that is the strict transform of $\alpha^*D$. 
We can assume that $\widetilde f:\widetilde X 
\longrightarrow Y'$ and $D'$ satisfy the conditions 
and the assumptions in 
\ref{setting} and Theorem \ref{can-ext} for a suitable 
simple normal crossing divisor $\Sigma'$ on $Y'$. 
By Proposition \ref{relative} and 
Corollary \ref{useful}, we can check that 
$R^i\widetilde f_*\omega_{\widetilde X}(D')\simeq R^i{f'}_*\omega_
{X'}(\alpha ^*D)$ for $i\geq 0$. 
We note that $(X',\alpha^*D)$ is lc 
and every center of log canonical singularities for 
the pair $(X',\alpha ^*D)$ is dominant onto $Y'$. 
Thus, $R^i{f'}_*\omega_
{X'}(\alpha ^*D)$ is locally free for $i\geq 0$. 
Since $R^if_*\omega_X(D)$ is a direct summand of $\tau_*
R^i{f'}_*\omega_{X'}(\alpha ^*D)$, we obtain 
that $R^if_*\omega_X(D)$ is locally free for $i\geq 0$. 
So, we complete the proof. 
\end{proof}

\subsection{Semi-positivity theorem}\label{subsec3.2}

In this subsection, we prove the semi-positivity of 
$R^if_*\omega_{X/Y}(D)$ on suitable assumptions. 
It is a generalization of 
Fujita-Kawamata's semi-positivity 
theorem and related to \cite[Theorem 32]{ka}. 
See Caution \ref{ka-ni-tuite} below. 

First, let us recall the definition of semi-positive vector bundles. 

\begin{defn}[Semi-positivity]\label{semi-def}
Let $V$ be a complete variety and $\mathcal E$ a 
locally free sheaf on $V$. 
We say that $\mathcal E$ is {\em{semi-positive}} 
if and only if the tautological line bundle 
$\mathcal O_{\mathbb P_V(\mathcal E)}(1)$ 
is nef on $\mathbb P_V (\mathcal E)$. 
We note that $\mathcal E$ is semi-positive if and 
only if for every complete curve $C$ and 
morphism $g:C\longrightarrow V$ 
every quotient line bundle of $g^*\mathcal E$ 
has non-negative degree. 
\end{defn}

We collect the basic properties of semi-positive vector bundles 
for the readers' convenience. We omit the proof here. Details are 
left to the readers. See the corresponding part of 
\cite{la}.   

\begin{prop}[Properties of semi-positive vector bundles]\label{propert}
Let $V$ be a complete variety. Then we have the 
following properties{\em{:}} 
\begin{itemize}
\item[(i)] Let $\mathcal E_1$ and $\mathcal E_2$ be 
vector bundles on $V$. 
Then the direct sum $\mathcal E_1\oplus \mathcal E_2$ is 
semi-positive if and only if both $\mathcal E_1$ and 
$\mathcal E_2$ are semi-positive. 
\item[(ii)] A vector bundle $\mathcal E$ on $V$ is 
semi-positive 
if and only if so is $S^k\mathcal E$ for every $k$, 
where $S^k\mathcal E$ is the $k$-th symmetric product of $\mathcal E$. 
\item[(iii)] If $\mathcal E$ is semi-positive and 
$\mathcal F$ is a semi-positive $($resp.~an ample$)$ vector bundle 
on $V$, then $\mathcal E\otimes \mathcal F$ is semi-positive $($resp.~
ample$)$. 
\item[(iv)] Any tensor product or exterior product of semi-positive 
vector bundles is semi-positive. 
\end{itemize}
\end{prop}

\begin{rem}\label{nef-ve}
It is obvious that a line bundle $\mathcal L$ 
is semi-positive if and only 
if $\mathcal L$ is nef. 
We note that, in some literatures (for example, \cite{la}), 
semi-positive vector bundles are called 
{\em{nef vector bundles}}. 
\end{rem}

The following lemma, which is not difficult to 
prove, will be used in the proof of Theorem 
\ref{semi}. We leave the details to the readers. 

\begin{lem}[Extension of semi-positive vector bundles]\label{ext}
Let $Y$ be a complete variety. 
Assume that there exists a short exact sequence on $Y${\em{:}}  
$$
0\longrightarrow \mathcal E'\longrightarrow 
\mathcal E\longrightarrow \mathcal E''\longrightarrow 0
$$ 
such that both $\mathcal E$ and $\mathcal E''$ are semi-positive 
vector bundles. 
Then so is $\mathcal E$. 
\end{lem}

The next theorem is the main theorem of this subsection. 
We will use it in Section \ref{sec4} to generalize Kawamata's 
positivity theorem. 

\begin{thm}[Semi-positivity theorem]\label{semi}
Let $f:X\longrightarrow Y$ be a projective surjective morphism 
between non-singular varieties over $\mathbb C$. 
Let $D$ be a simple normal crossing divisor on $X$ such 
that $D$ is strongly horizontal. 
Assume that there is a non-empty Zariski open set $Y_0$ 
of $Y$ such that $\Sigma:=Y\setminus Y_0$ is a simple 
normal crossing divisor on $Y$ 
and that $f_0:X_0\longrightarrow Y_0$ 
is smooth and $D_0$ is 
relatively normal crossing over 
$Y_0$, where $X_0:=f^{-1}(Y_0)$, 
$f_0:=f|_{X_0}$ and $D_0:=D\cap X_0$. 
Let $\omega_{X/Y}:=\omega_X\otimes f^*\omega^{-1}_Y$ and 
$d:=\dim f$. 
Assume that all the local monodromies of the 
local system $R^{d+i}f_{0*}\mathbb C_{X_0-D_0}$ 
around $\Sigma$ are unipotent. 
Let $Z$ be a complete subvariety of $Y$. 
Then the restriction $R^if_*\omega_{X/Y}(D)|_Z$ 
is semi-positive. 
In particular, if $Y$ is complete, then 
$R^if_*\omega_{X/Y}(D)$ is semi-positive. 
\end{thm}
\begin{proof} 
Let 
$$0\subset \cdots \subset W_k\subset 
W_{k+1}\subset \cdots \subset 
\mathcal H_{0}^{d+i}:=R^{d+i}f_{0*}\mathbb C_{X_0-D_0}
$$ 
be the weight filtration 
of the gradedly polarized variation of mixed Hodge structures and 
$$
0\subset \cdots \subset \widetilde 
W_k\subset \widetilde W_{k+1}\subset \cdots \subset 
\mathcal H_Y^{d+i}
$$ 
be the canonical extension 
of the above weight filtration. 
Then, the vector bundle $F^d(\mathcal H_Y^{d+i})$ induces the canonical 
extension of the bottom Hodge filtration on each $\Gr_k^{\widetilde W}
(\mathcal H_Y^{d+i})$. 
Therefore, Lemma \ref{ext}, \cite[Theorem 2]{ka2} 
and \cite[\S 4]{ka}
\footnote{Kawamata's proof of semi-positivity theorem 
heavily relies on the asymptotic behavior 
of the Hodge metric near a puncture. It is not so easy for the 
non-expert to take it out from \cite[\S 6]{schmid}. 
We recommend the readers to see \cite[Sections 2, 3]{p} or \cite
{zu}. Section 4 of \cite{f3} is an exposition of Fujita-Kawamata's 
semi-positivity theorem.} 
imply that 
$F^d(\mathcal H_Y^{d+i})|_Z$ is semi-positive. 
On the other hand, by Theorem \ref{can-ext}, 
we have the following isomorphism 
$$
R^if_*\omega_{X/Y}(D)\simeq F^d(\mathcal H_Y^{d+i}).   
$$ 
Thus, we obtain that $R^if_*\omega_{X/Y}(D)|_Z$ is semi-positive. 
\end{proof}

\begin{cau}\label{ka-ni-tuite}
The semi-positivity of $f_*\omega_{X/Y}(D)$ in Theorem 
\ref{semi} is very similar to \cite[Theorem 32]{ka}. 
Unfortunately, our theorem:~Theorem \ref{semi}, does not 
contain Theorem 32 in \cite{ka}. 
Note that we assume that $D$ is strongly horizontal, in 
particular, $D$ has no vertical components. 
This 
assumption is a little stronger than Kawamata's. 
\end{cau}

\section{A generalization of Kawamata's positivity theorem}\label{sec4}

\subsection{Positivity theorem}\label{subsec4.1}

The following theorem is one of the main results in this paper. 
It is a slight but important 
generalization of Kawamata's positivity theorem 
(see \cite[Theorem 2]{ka3}). See also \cite{katsuika} 
for the case where the fibers are curves. 

\begin{thm}[A generalization of Kawamata's positivity 
theorem]\label{kpt}
Let $f: X \longrightarrow B$ be a projective surjective morphism 
between non-singular quasi-projective varieties with
connected fibers.
Let $P = \sum P_j$ and $Q = \sum_{l} Q_{l}$ be
simple normal crossing divisors on $X$ and $B$, 
respectively, such that
$f^{-1}(Q) \subset P$ and 
$f$ is smooth over $B \setminus Q$.
Let $D = \sum_j d_jP_j$ be a $\mathbb Q$-divisor on $X$,
where $d_j$ may be positive, zero or negative, 
which satisfies the following conditions{\em{:}}
\begin{itemize}
\item[(1)] 
$f: \Supp(D^h) \longrightarrow B$ is relatively normal 
crossing over $B \setminus Q$, and
$f(\Supp(D^v)) \subset Q$.  

\item[(2)] $d_j \leq 1$ if $P_j$ is horizontal. 

\item[(3)] $\dim_{k(\eta)}f_*\mathcal O_X
(\ulcorner - D^{<1} \urcorner)\otimes 
_{\mathcal O_B}k(\eta)=1$ for the generic point $\eta$ of $B$. 

\item[(4)] $K_X + D \sim_{\mathbb Q} f^*(K_B + L)$ for some 
$\mathbb Q$-divisor $L$ on B.
\end{itemize} 
We put 
$$
\gamma_l:=\max \{\gamma\in \mathbb R | (X, D+\gamma f^*Q_l) 
{\text{is sub lc over the generic point of}} \ Q_l\},   
$$ 
that is, $\gamma _l$ is the log canonical threshold for the 
pair $(X,D)$ with respect to $f^*Q_l$ over the 
generic point of $Q_l$. 
We put 
$$\delta _l:=1-\gamma _l$$ 
and define 
\begin{eqnarray*}
\Delta_0:=\sum _l\delta_l Q_l \\
M:=L-\Delta_0. 
\end{eqnarray*}
Then $M\cdot C\geq 0$ for every projective curve 
$C$ on $B$. In particular, 
if $B$ is projective, then $M$ is nef in the usual sense. 
\end{thm}

\begin{rem}\label{kawamata}
We use the same notation as in Theorem \ref{kpt}. 
Let 
\begin{eqnarray*}
f^{*}Q_{l} &=&\sum_j w_{l j}P_j, \\
\bar d_j &=& \frac {d_j + w_{l j} - 1}{w_{l j}}
\text{ if } f(P_j) = Q_{l}. 
\end{eqnarray*} 
Then we can check easily that 
$$
\delta_{l} =\max\{ \bar d_j\ |\ f(P_j) = Q_{l}\}   
$$ by the definition of the log canonical threshold. 
So, our definition of $M$ coincides with Kawamata's. 
\end{rem}

\begin{rem}\label{hitsuyou}
In \cite[Theorem 2]{ka3}, it is assumed that $d_j< 1$ 
for all $j$. On this assumption, it is obvious that 
$D^{<1}=D$. 
So, our theorem contains original Kawamata's positivity 
theorem. 
\end{rem}

Before we give the proof of Theorem \ref{kpt}, 
we recall the following well-known lemma. 
It may help the readers to understand the 
proof of Theorem \ref{kpt}, which is related to 
the toroidal geometry. 
We learned it from \cite{ak}. 
Proposition 3.1 and Lemma 6.2 of \cite{ak} 
are useful for us. 

\begin{lem}\label{toric}
Let $X$ be a toric variety and $D$ the complement of 
the big torus. 
Then $(X,D)$ is lc. 
\end{lem}

The following proof is essentially the same as Kawamata's. 
We repeat his arguments in details for the 
readers' convenience. 
In his proof, he used Fujita-Kawamata's semi-positivity 
theorem. 
On the other hand, we apply Theorem \ref{semi}. 

\begin{proof}[Proof of {\em{Theorem \ref{kpt}}}] 
By replacing $D$ by $D - f^*\Delta_0$, 
we can assume that $\Delta_0 = 0$.
Then we have an inequality $d_j \le 1 - w_{l j}$ for $f(P_j) = Q_{l}$, 
and the equality holds for some $j$ for each $l$.

By the semi-stable reduction theorem and Kawamata's 
covering trick, 
we obtain a semi-stable reduction 
in codimension one in the following sense:
~there exists a finite morphism 
$h: B' \longrightarrow B$ from a non-singular quasi-projective 
variety $B'$ with a simple normal crossing divisor 
$Q':= \Supp (h^*Q) = \sum_{l'} Q'_{l'}$ 
such that the induced morphism 
$f': X' \longrightarrow B'$ from a desingularization $X'$ 
of $X \times_B B'$ is semi-stable over the generic points of $Q'$.
Let $g: X' \longrightarrow X$ be the induced morphism.
We may assume that $P' = \text{Supp}(g^*P) = \sum_{j'} P'_{j'}$ 
is a simple normal crossing divisor again:
$$
\begin{CD}
X     @<g<< X' \\
@VfVV       @VV{f'}V \\
B     @<h<< B'.
\end{CD}
$$

Let $Z \subset B'$ be a closed subset of codimension two 
or larger that is 
contained in $Q'$ and such that $Q' \setminus Z$ is smooth and
$f'$ is semi-stable over $B' \setminus Z$.
We can define naturally a 
$\mathbb Q$-divisor $D' = \sum_{j'}d'_{j'}P'_{j'}$ on 
$X'$ such that 
$K_{X'} + D' \sim_{\mathbb Q} f^{\prime *}(K_{B'} + h^*L)$.
We calculate the coefficients $d'_{j'}$.
If $P'_{j'}$ is horizontal and $g(P'_{j'}) = P_j$, then $d'_{j'} = d_j$.
We are not concerned with those $P'_{j'}$ such that $f'(P'_{j'}) \subset Z$.

We consider the case where $f'(P'_{j'}) = Q'_{l'}$. 
Then, we have an inequality $d'_{j'} \le 0$, 
and the equality holds for some $j'$ for each $l'$.
The inequality $d'_{j'} \le 0$ holds for any $j'$ such that $f'(P'_{j'}) = 
Q'_{l'}$, 
because it is stable under the blow-ups of $X$. 
We note that $\gamma _l=1-\delta_l$ is the log canonical 
threshold for the pair $(X,D)$ with respect to 
the divisor $f^*Q_l$ over the generic point of $Q_l$. 

Let $m$ be a positive integer such that $mL$ is a Weil divisor, 
$X'_m := X' \times_{B'} \cdots \times_{B'} X'$ the  
$m$-tuple fiber product of $X'$ over $B'$, and
$f'_m: X'_m \longrightarrow B'$ the projection.
Then 
$X'_m$ has only Gorenstein toric singularities over $B' \setminus Z$.
Note that $f_m'$ is weakly semi-stable over $B'\setminus Z$ by 
\cite[Definition 0.1, Corollary 1.6, Lemma 6.2]{ak}. 
Let $D'_m = \sum_{i = 1}^m p_i^*D'$ be the sum 
over the $i$-th projections $p_i: X'_m \longrightarrow X'$.
Then we have $K_{X'_m} + D'_m \sim_{\mathbb Q} f_m^{\prime *}(K_{B'} + mh^*L)$ 
over $B' \setminus Z$. 

Let $r$ be the smallest positive integer such that 
$r(K_{X'_m} + D'_m - f_m^{\prime *}(K_{B'} + mh^*L)) \sim 0$ 
over $B' \setminus Z$, and $\theta$ a 
non-zero rational function such that 
$$
\xdiv(\theta)=r(K_{X'_m} + D'_m - f_m^{\prime *}(K_{B'} + mh^*L)) 
$$ over $B' \setminus Z$. 
Let $\pi: 
\widetilde X \longrightarrow X'_m$ be 
the normalization of the main irreducible 
component $X^{\prime o}_m$ of $X'_m$ in the field 
$\mathbb C(X^{\prime o}_m)(\root{r}\of{\theta})$, and 
$\widetilde f: \widetilde X \longrightarrow B'$ the induced morphism.
We note that 
$\pi$ may ramify along the support of $D'_m$ over $B' \setminus Z$,
and $\widetilde f$ may have non-connected fibers. 
Then we obtain 
$$
\pi_*\mathcal O_{\widetilde X}\simeq 
\bigoplus^{r-1}_{k=0}\mathcal O_{X_m'}(kK_{X_m'}+
\llcorner kD_m'\lrcorner
-k f_m^{\prime *}(K_{B'} + mh^*L))  
$$ 
over $B'\setminus Z$. By the duality, 
$$
\pi_*\omega_{\widetilde X}\simeq 
\bigoplus^{r-1}_{k=0}\mathcal O_{X_m'}
((1-k)K_{X_m'}+\ulcorner -kD_m'\urcorner
+k f_m^{\prime *}(K_{B'} + mh^*L))  
$$ 
over $B'\setminus Z$. 
Let $\rho:W\longrightarrow \widetilde X$ 
be a resolution of singularities 
such that 
$\lambda:=f'_m\circ \pi\circ \rho:W \longrightarrow B'$ 
is smooth over $B'\setminus Q'$. 
We note that since $\widetilde X\longrightarrow B'$ is 
equisingular over $B'\setminus Q'$ by the construction, 
there is a simultaneous resolution over 
$B'\setminus Q'$ by the canonical desingularization 
theorem. 
We can further assume that there exists a simple normal crossing 
divisor $F$ on $W$ such that $F$ is relatively normal crossing over 
$B'\setminus Q'$, strongly horizontal over $B'$, and 
$\rho_*\omega_W(F)\simeq \omega_{\widetilde X}(\pi^*E)$ over 
$B'\setminus Z'$, where $E:=D_m'-{D_m'}^{<1}$. 
We note that $(\widetilde X, \pi^*E)$ is log canonical 
such that all the discrepancies of the divisors whose centers are not 
dominant onto $B'$ are non-negative
over $B'\setminus Z$. 

We consider the local system $R^d{\lambda_0}_*\mathbb C_{W_0-F_0}$, 
where $B'_0:=B'\setminus Q'$, $W_0:=\lambda^{-1}(B'_0)$, 
$\lambda_0:=\lambda|_{W_0}$, $F_0:=F\cap W_0$, and 
$d:=\dim \lambda$. 
By the covering trick again, 
we construct a finite morphism $h': B'' \longrightarrow B'$ from a 
non-singular quasi-projective variety $B''$
with a simple normal crossing divisor 
$Q'' := \text{Supp}(h^{\prime *}Q') = \sum_{l''} Q''_{l''}$ 
such that the induced morphism $\lambda ': 
W' \longrightarrow B''$ 
from a desingularization $\mu: W' \longrightarrow W \times_{B'} B''$ 
has unipotent local monodromies on 
the local system $(h')^{-1}R^d{\lambda _0}_*\mathbb C_{W_0-F_0}$ 
around the irreducible components of
$Q''$. 
It is a unipotent reduction with respect to 
the local system $R^d{\lambda_0}_*\mathbb C_{W_0-F_0}$. 

Let $X''$ be the normalization of the main irreducible components of the
fiber product $X'_m \times_{B'} B''$.
Let $f'': X'' \longrightarrow B''$ 
and $g': X'' \longrightarrow X'_m$ be the induced morphisms.
Since $f'$ is semi-stable over $B' \setminus Z$, we have
$g^{\prime *}K_{X'_m/B'} = K_{X''/B''}$ over $B'' \setminus Z'$ for
$Z' = h^{\prime -1}(Z)$.
Thus $\theta' = g^{\prime *}\theta$ is a rational function such that 
$$
\xdiv(\theta')=r(K_{X''/B''} + {g'}^{*}D'_m 
- f^{'' *}mh^{'*}h^*L) 
$$ 
over $B'' \setminus Z'$. 
Let $\widetilde X'$ be the normalization of $X''$ 
in the field $\mathbb C(X'')(\root{r}\of{\theta'})$.
We note that $f''$ is weakly semi-stable over $B''\setminus Z'$ by 
\cite[Lemma 6.2]{ak}. 
Let $\pi': \widetilde X' \longrightarrow X''$ and 
$\widetilde f': \widetilde X' \longrightarrow B''$ be the induced morphisms:
$$
\CD
W @<<< W'\\ 
@V{\rho}VV  @VV{\rho'}V \\
\widetilde X   @<{\widetilde g}<< \widetilde X' \\
@V{\pi}VV                 @VV{\pi'}V \\
X'_m       @<{g'}<<       X'' \\
@V{f'_m}VV                @VV{f''}V \\
B'         @<{h'}<<       B''.
\endCD
$$ 
Since $\pi'$ ramifies only along the support of $D'' = g^{\prime *}D'_m$ over 
$B'' \setminus Z'$,
$\widetilde X'$ has only toric singularities there. 
We have 
$$
\pi'_*\omega_{\widetilde X'}\simeq 
\bigoplus^{r-1}_{k=0}\mathcal O_{X''}((1-k)K_{X''}+\ulcorner -kD''
\urcorner
+k {f''}^{*}(K_{B''} + m{h'}^{*}h^*L))  
$$ 
over $B''\setminus Z'$ by 
the same arguments as above. 
We put $E':=D''-{D''}^{<1}$. 
We note that 
$(\widetilde X',{\pi'}^*E')$ is lc such that 
all the discrepancies of the divisors whose centers are not 
dominant onto $B''$ are non-negative over $B''\setminus
Z'$. So, we can assume that there exists a simple normal 
crossing divisor $F'$ on $W'$ such that 
$F'$ is strongly horizontal, relatively normal crossing over
$B''\setminus Q''$, and that $\rho' _*\omega_{W'}(F')
\simeq \omega_{\widetilde X'}({\pi'} ^*E')$. 
We note that since $(\widetilde X',{\pi'}^*E')$ 
is equisingular over 
$B''\setminus Q'$, there is a simultaneous resolution over 
$B''\setminus Q'$ by the canonical desingularization theorem. 
We can assume that $W'=W\times _{B'}B''$ and 
$F'=F\times_{B'}B''$ over $B''\setminus Q''$. 
We put $\lambda'=f''\circ \pi'\circ\rho'$, 
$W'_0:=({\lambda'})^{-1}(B''\setminus Q'')$, 
$\lambda'_0:=\lambda'|_{W'_0}$, and 
$F'_0:=F'\cap W'_0$. 
Then, $\lambda'_*\omega_{W'}(F')
\simeq {f''}_*\pi'_*\omega_{\widetilde X'}({\pi'}^*E')$ 
is semi-positive 
when it is restricted to 
a projective subvariety 
by Theorem \ref{semi} 
since all the local monodromies on 
$R^d{\lambda'_0}_*\mathbb C_{W'_0-F'_0}$ 
around $Q''$ are unipotent. 

By the above argument, 
we have that 
$$
f''_*\mathcal O_{X''}(\ulcorner -{D''}^{<1}\urcorner)
\otimes \mathcal O_{B''}(m{h'}^{*}h^*L)
$$ 
is a direct summand of 
$\lambda'_*\omega_{W'}(F')$. 
Since we have $d'_{j'}\leq 0$ with the equality for some 
$j'$ for each $l'$, we have 
$f''_*\mathcal O_{X''}(\ulcorner -{D''}^{<1}\urcorner)\simeq 
\mathcal O_{B''}$. 
So, 
$\mathcal O_{B''}(m{h'}^{*}h^*L)$ is a direct summand of 
$\lambda'_*\omega_{W'}(F')$. 
Thus we obtain that $M\cdot C=L\cdot C\geq 0$ 
for every projective curve $C$ on $B$. 
\end{proof}

\subsection{Applications}\label{subsec4.2}

The following is a slight generalization of \cite[Theorem 0.2]{f1}. 
We explain the proof in details for the readers' convenience. 

\begin{thm}\label{mt}
Let $(X,\Delta)$ be a proper sub lc pair 
and $f:X\longrightarrow S$ a proper surjective morphism 
onto a normal variety $S$ with connected fibers 
such that $\Delta$ is strongly horizontal with respect to $f$. 
Assume that 
$\dim _{k(\eta)} f_{*}{\mathcal O}_X (\ulcorner -\Delta^{<1}\urcorner)
\otimes_{{\mathcal O}_{S}} k(\eta)=1$, 
where $\eta$ is the generic point of $S$. 
We further 
assume that there exists a 
$\bQ$-Cartier $\bQ$-divisor $A$ 
on $S$ such that $K_X+\Delta\sim_{\bQ}f^{*} A$. 
Let $H$ be an ample Cartier divisor on $S$, and $\epsilon$ a 
positive rational number. 
Then there exists a $\bQ$-divisor $B$ on $S$ 
such that 
$$
K_S+B\sim_{\bQ}A+\epsilon H ,
$$
$$
K_X+\Delta+\epsilon f^{*} H\sim _{\bQ} f^{*} (K_S+B),
$$ 
and that the pair $(S,B)$ is sub klt. 

Furthermore, if 
$f_{*}{\mathcal O} _{X} (\ulcorner -\Delta^{<1}
\urcorner) \simeq {\mathcal O}
_S$, 
then we can make $B$ effective, that is, $(S,B)$ is klt. 
In particular, $S$ has only rational singularities. 
\end{thm}

\begin{proof} (cf.~Proof of \cite[Theorem 2]{n2}). 
By using the desingularization theorem, 
we have the following commutative diagram:
$$
\begin{CD}
Y @>\text{$\nu$}>> X \\
  @V\text {$g$}VV  @VV\text {$f$}V  \\
T @>\text{$\mu$}>>S, 
\end{CD}
$$
where 
\begin{enumerate}
\item [(i)] $Y$ and $T$ are non-singular projective varieties, 
\item [(ii)] $\nu$ and $\mu$ are projective birational morphisms, 
\item [(iii)] we define $\bQ$-divisors $D$ and $L$ on $Y$ and $T$ 
by the following relations:
$$
K_Y+D=\nu^{*}(K_X+\Delta),
$$
$$
K_T+L  \sim _{\bQ} \mu^{*} A,
$$
\item [(iv)] there are simple normal crossing divisors $P$ and $Q$ on 
$Y$ and $T$ such that they satisfy the conditions 
of Theorem \ref {kpt} and 
there exists a set of positive rational numbers 
$\{s_l\}$ such that $\mu^{*} H-\sum_{l} s_l Q_l$ is ample. 
\end{enumerate}
By the construction, the conditions (1) and (4) of Theorem 
\ref {kpt} are satisfied. 
Since $(X,\Delta)$ is sub lc, the condition (2) of Theorem 
\ref {kpt} is satisfied. 
The condition (3) of Theorem 
\ref {kpt} can be checked by the following claim. 
Note that $\mu$ is birational. 
We put $h:=f\circ \nu$.
\begin{claim}[A]\label{A}
${\mathcal O}_{S} \subset h_{*}{\mathcal O}_{Y}(\ulcorner -D^{<1}
\urcorner) 
\subset f_{*}{\mathcal O}_{X}(\ulcorner -\Delta^{<1}\urcorner)$.
\end{claim}

\begin{proof}[Proof of {\em{Claim (A)}}]
First, we have  
${\mathcal O} _{S} 
\simeq h_{*}{\mathcal O} _{Y}
\subset h_{*}{\mathcal O}_{Y}(\ulcorner -D^{<1}\urcorner)$, 
since ${\mathcal O}_{Y}\subset {\mathcal O}_{Y}(\ulcorner -D
^{<1}\urcorner)$. Next, we have 
$$
\begin{array}{clcl}
&\Gamma(U, \nu_{*}{\mathcal O}_{Y}(\ulcorner -D^{<1}\urcorner))& 
\subset & \Gamma(U \setminus 
Z, \nu_{*}{\mathcal O}_{Y}(\ulcorner -D^{<1}\urcorner))\\
=& \Gamma(U \setminus
Z, {\mathcal O}_{X}(\ulcorner -\Delta^{<1}\urcorner))
&=&\Gamma(U, {\mathcal O}_{X}(\ulcorner -\Delta^{<1}\urcorner)), 
\end{array}
$$
where $U$ is a Zariski open set of $X$ and $Z:=\nu (\Exc (\nu))$. 
So we have $\nu_{*}{\mathcal O}_{Y}(\ulcorner -D^{<1}\urcorner)
\subset {\mathcal O}
_{X}(\ulcorner -\Delta^{<1}\urcorner)$. 
Then we obtain $h_{*}{\mathcal O}_{Y}(\ulcorner -D^{<1}\urcorner) 
\subset f_{*}{\mathcal O}_{X}(\ulcorner -\Delta^{<1}\urcorner)$. 
We complete the proof of Claim (A). 
\end{proof}

So we can apply Theorem \ref {kpt} to $g:Y\longrightarrow T$. 
The divisors $\Delta_{0}$ and $M$ are as in Theorem \ref {kpt}. 
Then $M$ is nef. 
Since $M$ is nef, we have that 
$$
M+\epsilon\mu^{*}H-\epsilon'\sum_{l} s_l Q_l 
$$ 
is ample for $0<\epsilon'\leq \epsilon$. 
We take a general Cartier divisor 
$$ 
F_0\in |m(M+\epsilon\mu^{*}H-\epsilon'\sum_{l} s_l Q_l )|
$$
for a sufficiently large and divisible integer $m$. 
We can assume that $\Supp (F_0 \cup \sum _l Q_l)$ is a simple normal 
crossing divisor. 
And we define $F:=(1/m)F_{0}$. 
Then 
$$
L+\epsilon\mu^{*}H\sim_{\bQ} F+\Delta_0 + \epsilon'\sum_{l} s_l Q_l .
$$
Let $B_0:= F+\Delta_0 + \epsilon'\sum_{l} s_l Q_l$ and 
$\mu_* B_0=B$. 
We have $K_T+B_{0}=\mu ^{*} (K_S+B)$. 
By the definition of $\Delta_0$ and 
the assumption that $\Delta$ is 
strongly horizontal with respect to $f$, 
$\llcorner\Delta_0\lrcorner\leq 0 $. 
So $\llcorner F+\Delta_0 + \epsilon'\sum_{l} s_l Q_l \lrcorner \leq 0$ 
when $\epsilon'$ is small enough. 
Then $(S,B)$ is sub klt. 
By the construction we have  
$$
K_S+B\sim_{\bQ}A+\epsilon H ,
$$
$$
K_X+\Delta+\epsilon f^{*} H\sim_{\bQ} f^{*}(K_S+B).
$$
If we assume furthermore that 
$f_{*}{\mathcal O}_X (\ulcorner -\Delta^{<1}
\urcorner)\simeq {\mathcal O}_S$, 
we can prove the following claim. 
For the notation:~$\delta_l, w_{lj}$, and $\bar d_j$, 
see Theorem \ref{kpt}. 

\begin{claim}[B]
If $\mu_* Q_l\ne 0$, then $\delta_l\geq 0$.
\end{claim}
\begin{proof}[Proof of {\em{Claim (B)}}]
If $\ulcorner -d_{j}\urcorner\geq w_{lj}$ for every $j$, 
then $\ulcorner -D^{<1}\urcorner\geq g^{*}Q_{l}$. 
So $g_{*} {\mathcal O} _{Y}( \ulcorner -D^{<1}\urcorner) \supset 
{\mathcal O} _{T} (Q_l)$. 
Then $\mathcal O _{S} \simeq h_{*}{\mathcal O} _{Y}( \ulcorner -D^{<1}
\urcorner)
\supset \mu _{*} {\mathcal O} _{T} (Q_l)$ by Claim (A). 
It is a contradiction. 
So we have that $\ulcorner -d_j\urcorner <w_{lj}$ for some $j$. 
Since $w_{lj}$ is an integer, we have that 
$-d_j+1\leq w_{lj}$. 
Then $\bar d_j\geq 0$. We get $\delta_l\geq 0$. 
\end{proof}

Therefore, $B$ is effective if  
$f_{*}{\mathcal O}_X (\ulcorner -\Delta^{<1}\urcorner)\simeq 
{\mathcal O}_S$. 
This completes the proof. 
\end{proof}

\begin{rem}
Under the same notation and assumptions as in Theorem \ref{mt}, 
we further assume that $S$ is projective and 
$f_{*}{\mathcal O} _{X} (\ulcorner -\Delta^{<1}
\urcorner) \simeq {\mathcal O}
_S$. Then the (generalized) cone theorem holds 
on $S$ with respect to $A$. 
For the details, see \cite[Section 4]{f1}. 
We note that \cite[Section 4]{f1} was completely 
generalized in \cite{am}. In his notation, 
$(S,A)$ is a projective {\em{quasi-log variety}}. 
\end{rem}

\section{Log canonical bundle formula}\label{sec5}

In \cite[Section 4]{fm}, we formulated and proved 
a log canonical bundle formula for klt pairs. 
In this section, we give a log canonical bundle formula 
for log canonical pairs. The main theorem of 
this section is Theorem \ref{mainthm}. 

\begin{say5}\label{(s4.1)}
Let $f : X \longrightarrow S$ be
a proper surjective morphism of a normal variety $X$
of dimension $n=m+l$ to a non-singular $l$-fold
$S$ such that
\begin{enumerate} 
\item[(i)] $(X,\Delta)$ is a sub lc pair
(assumed lc from \ref{(ss4.1)} and on), and

\item[(ii)] the generic fiber $F$ of $f$ is a
geometrically irreducible variety
with $\kappa(F, (K_X+\Delta)|_F) = 0$.
We fix the smallest positive integer $b$
such that the $f_{*}{\mathcal O}_{X}(b(K_X+\Delta))\ne 0$.
\end{enumerate}
\end{say5}

\begin{prop5}[{\cite[Propositions 2.2, 4.2]{fm}}]\label{(x2.2)} 
There exists one and only one
$\bQ$-divisor $D$ 
modulo linear equivalence on $S$ with a
graded $\xO_S$-algebra isomorphism
$$
\bigoplus_{i \ge 0} \xO_S(\llcorner iD \lrcorner)
\cong \bigoplus_{i \ge 0} 
(f_*\xO_X(\llcorner ib(K_X+\Delta)\lrcorner-ibf^*K_S))^{**}.
$$
Furthermore, the above isomorphism induces the equality{\em{:}} 
$$
b(K_X+\Delta) = f^*(b K_S+D)+B^{\Delta},
$$
where $B^{\Delta}$ is a $\bQ$-divisor on $X$ such that 
$f_*\xO_X(\llcorner iB_+^{\Delta}\lrcorner)\simeq \xO_S$ for 
every $i>0$ and $\codim_S(f(\Supp B_-^{\Delta})) \ge 2$.
\end{prop5}

\begin{rem5}\label{iruka}
In Proposition \ref{(x2.2)}, we note that 
for an arbitrary open set $U \subset S$, $D|_{U}$ and 
$B^{\Delta}|_{f^{-1}(U)}$ depend only on $f|_{f^{-1}(U)}$
and $\Delta|_{f^{-1}(U)}$.
\end{rem5}

\begin{defn5}[$L_{X/Y}^{\log}$]\label{(x2.3)}
We denote the $\mathbb Q$-divisor 
$D$ given in Proposition \ref{(x2.2)} by $L_{(X,\Delta)/S}$ or
simply by $L_{X/S}^{\log}$ if there is no danger of confusion.
\end{defn5}

\begin{defn5}[$s_P$ and $t_P$]\label{(x2.3.1)}
Let $P$ be a prime divisor on $S$. 
We set $s_P^{\Delta}:=b(1-t_P^{\Delta})$, where
$t_P^{\Delta}$ is the log canonical threshold of ${f}^{*}P$ with
respect to $(X, \Delta-(1/b)B^{\Delta})$ over the
generic point $\eta_P$ of $P$:
$$
t_P^{\Delta}:=\max\{t\in \mathbb R\ |\
(X, \Delta-(1/b)B^{\Delta}+tf^{*}P)
\text{ is sub lc over } \eta_P\}.
$$
\end{defn5}
Note that $s_P\in \mathbb Q$ and 
$s_P^{\Delta}\ne 0$ only for 
a finite number of codimension one points $P$
because there exists a nonempty Zariski open set $U\subset S$ 
such that $s_P^{\Delta}=0$ for every prime divisor $P$ with 
$P\cap U\ne \emptyset$. We may simply write $s_P$ 
rather than $s_P^{\Delta}$ if there is
no danger of confusion.
We note that $s_P^{\Delta}$ depends only on
$f|_{f^{-1}(U)}$ and $\Delta|_{f^{-1}(U)}$ where
$U$ is an open set containing $P$.

\begin{defn5}[Log-semistable part $L_{X/S}^{\log,ss}$]\label{(x2.3.2)}
We set 
$$L_{(X,\Delta)/S}^{ss}:=L_{(X,\Delta)/S}-\sum_P 
s_P^{\Delta} P
$$
and call it the {\it log-semistable part} of $f$.
We may simply denote it by $L_{X/S}^{\log,ss}$ if
there is no danger of confusion.
\end{defn5}
\begin{rem5}
We note that $D$, $L_{(X,\Delta)/S}$, $s_P^{\Delta}$,
$t_P^{\Delta}$ and $L_{(X,\Delta)/S}^{ss}$
are birational invariants of $(X,\Delta)$ over $S$ in the
following sense. 
Let $(X',\Delta')$ be a projective sub lc pair
and  $\sigma: X' \longrightarrow X$ a birational morphism
such that $K_{X'}+\Delta' - \sigma^*(K_X+\Delta)$
is an effective $\sigma$-exceptional $\bQ$-divisor.
Then the above invariants for
$f\circ \sigma$ and $(X',\Delta')$ 
are equal to those for $f$ and $(X,\Delta)$.
\end{rem5}

\begin{say5}\label{58}
Putting the above symbols together, we have 
{\it the log canonical bundle formula}
for $(X,\Delta)$ over $S$:
\begin{equation}\label{lcbformula}
b(K_X+\Delta) = f^*(bK_S+L_{X/S}^{\log,ss})+
\sum_P s_P^{\Delta} f^*P
+ B^{\Delta},
\end{equation}
where $B^{\Delta}$ is a $\bQ$-divisor on $X$ such that 
$f_*\xO_X(\llcorner iB_+^{\Delta}\lrcorner)\simeq \xO_S\ 
(\forall i>0)$ and
$\codim_S(f(\Supp B_-^{\Delta})) \ge 2$.
\end{say5}

We need to pass to a certain 
birational model
$f' : X' \longrightarrow S'$ to understand the log-semistable part more
clearly and to make the log canonical bundle formula more useful.

\begin{say5}\label{(ss4.1)}
From now on, we assume that $(X,\Delta)$ is lc. 

By Proposition \ref{(x2.2)}, we have 
$$
K_X+\Delta -\frac{1}{b}B^{\Delta}
=f^*(K_S+\frac{1}{b}L^{\log}_{X/S}). 
$$ 
Let $g:Y\longrightarrow X$ be a log resolution of 
$(X,\Delta -(1/b)B^\Delta)$ 
with $G$ a $\bQ$-divisor on $Y$ such that 
$$
K_Y+G=g^*(K_X+\Delta-\frac{1}{b}B^{\Delta}). 
$$ 
We put $K_Y+\Theta=g^*(K_X+\Delta)$. 
Let $\Sigma \subset S$
be an effective divisor 
satisfying the following conditions; 
\begin{enumerate}
\item $h:=f\circ g$ is projective,  
\item $h$ is smooth and
$\Supp G^h$ is  relatively normal crossing over $S
\setminus \Sigma$, 
\item $h(\Supp G^v)\subset \Sigma$, and 
\item $f$ is flat over $S\setminus \Sigma$. 
\end{enumerate}
 
Let $\pi:S' \longrightarrow S$ be a proper birational morphism
from a non-singular quasi-projective variety 
such that
\begin{enumerate}
\item[(i)] $\Sigma':=\pi^{-1}(\Sigma)$ is a simple normal 
crossing divisor, 

\item[(ii)] $\pi$ induces an isomorphism
$S'\setminus \Sigma' \cong S\setminus \Sigma$, and

\item[(iii)] the irreducible component $X_1$ of $X \times_S S'$
dominating $S'$ is flat over $S'$.
\end{enumerate}
Let $X'$
be the normalization of
$X_1$, and $f':X'\longrightarrow S'$ 
the induced morphism.
Let $g':Y'\longrightarrow X'$ be a log resolution such that 
$Y'\setminus {h'}^{-1}({\Sigma'})\cong Y\times_{X}{X_1}
\setminus {\alpha}^{-1}({\Sigma'})$, 
where $h':=f'\circ g'$ and $\alpha:Y\times_{X}{X_1}\longrightarrow S'$. 
Let $\Theta'$ be the $\bQ$-divisor
on $Y'$ such that $K_{Y'}+\Theta'=(\tau\circ
g')^*(K_X+\Delta)$, where $\tau:X'\longrightarrow X$ is the induced
morphism.  
We put 
$$
K_{Y'}+G'=(\tau\circ g)^*(K_X+\Delta -\frac{1}{b}B^{\Delta}). 
$$ 
Furthermore, we can assume that 
$\Supp({h'}^{-1}(\Sigma')\cup G')$ is 
a simple normal crossing 
divisor, and $h'(\Supp{G'}^{v})\subset \Sigma'$. 
We note that 
$\Supp {G'}^h$ is relatively normal crossing 
over $S'\setminus \Sigma'$ by the construction. 

Later we treat horizontal or vertical divisors on $X, X'$
or $Y'$ over $S$ without referring to $S$.
Note that a $\bQ$-divisor on $X'$ or $Y'$
is horizontal (resp.~vertical) over $S$ if and only if it is 
horizontal (resp.~vertical) over $S'$.

We note that the horizontal part $(\Theta')^h_-$ of 
the negative part
$\Theta'_-$  of $\Theta'$ is 
$g'$-exceptional.  
$$
\begin{CD}
Y@<<<Y\times_{X}{X_1}@<<<Y'\\
@V{g}VV @VVV @VV{g'}V\\
X@<<< {X_1}    @<<<       {X'}\\
@V{f}VV @VVV   
@VV\text{$f'$}V\\ 
S    @<< \text{$\pi$}<   {S'} @= {S'}
\end{CD}
$$
\end{say5}

\begin{rem5}\label{misutta} 
The definition of $g$ in the above \ref{(ss4.1)} is 
slightly different from that in \cite{fm}. 
In \cite[4.4]{fm}, $g:Y\longrightarrow X$ is a 
log resolution of $(X,\Delta)$. 
However, it is better to assume that $g$ is 
a log resolution of $(X,\Delta-(1/b)B^\Delta)$ 
for the proof of Theorem \ref{nef}. 
See the conditions in Theorem \ref{kpt}. 
\end{rem5}

The following formula is the main theorem of this section. 
It is a slight generalization of \cite[Theorem 4.5]{fm}. 

\begin{thm5}[Log canonical bundle formula]\label{mainthm}
Under the above notation and assumptions, 
let $\Xi$ be a $\bQ$-divisor on $Y'$ such that 
$(Y',\Xi)$ is sub lc and
$\Xi-\Theta'$ is effective and exceptional over $X$. 
$($Note that $\Xi$ exists since $(X,\Delta)$ is lc.$)$ 
Then the log canonical bundle formula 
$$
b(K_{Y'}+\Xi) = (h')^*(bK_{S'}+L_{(Y',\Xi)/S'}^{ss})+\sum_P
s_P^{\Xi} (h')^*(P)+B^{\Xi}
$$
for $(Y',\Xi)$ over $S'$ has the following properties:
\begin{enumerate}
\item[(i)] $h'_*\xO_{Y'}(\llcorner iB^{\Xi}_+\lrcorner)
\simeq \xO_{S'}$ for all
$i>0$,
\item[(ii)] $B^{\Xi}_-$ is $g'$-exceptional and
$\codim_{S'} (h'(\Supp B^{\Xi}_-)) \ge 2$,
\item[(iii)] the following holds for every $i>0${\em{:}}
\begin{eqnarray*}
H^0(X, \llcorner ib(K_X+\Delta)\lrcorner)
&\simeq &H^0(Y', \llcorner ib(K_{Y'}+\Xi)\lrcorner)\\
&\simeq &H^0(S', \llcorner
ibK_{S'}+iL_{Y'/S'}^{\log,ss}+
\sum i \, s_P^{\Xi} P\lrcorner),
\end{eqnarray*}
\item[(iv)] $L_{Y'/S'}^{\log,ss}\cdot C\geq 0$ for 
every projective curve $C$ on $S'$, in particular, 
$L_{Y'/S'}^{\log,ss}$ is nef if $S$ is complete, and
\item[(v)] under the assumption that $(Y',\Xi)$ is lc, 
let $N$ be a positive integer such that $N
{h'}^*(L_{Y'/S'}^{\log,ss})$ and 
$bN\Xi^v$ are Weil divisors. Then
for each $P$, there exist
$u_P\in \xN$ and $v_P\in \mathbb Z_{\geq 0}$ such that
$0 \leq v_P \le bN$ and 
$$s_P = \frac{bNu_P-v_P}{Nu_P}.$$
\end{enumerate}
\end{thm5}

Before we give the proof, we note the following remark, 
which is obvious by the definition of $s_P$ and $t_P$. 

\begin{rem5}\label{tyotto}
Let $(Y,\Theta')$ be as in Theorem \ref{mainthm}. 
If $\Theta'$ is strongly horizontal with respect to 
$h'$, then $s_P<1$, equivalently, $v_P>0$. 
We note that $\Theta'$ is strongly horizontal 
with respect to $h'$ if and only if 
so is $\Delta$ with respect to $f$. 
\end{rem5}

\begin{proof}[Proof of {\em{Theorem \ref{mainthm}}}] 
First, (i) is obvious by the formula (\ref{lcbformula}) 
in \ref{58}. 
Similarly, (ii) follows because $(g')_*(B^{\Xi}_-)=0$ by the
equidimensionality of $f'$.

By (ii) and the conditions on $\Xi$, the following holds for all
$i>0$:
\begin{eqnarray*}
H^0(X, \llcorner ib(K_X+\Delta)\lrcorner)
&\simeq &H^0(Y', \llcorner ib(K_{Y'}+\Xi)\lrcorner)\\
&\simeq &H^0(Y', \llcorner ib(K_{Y'}+\Xi)+i B^{\Xi}_-\lrcorner).
\end{eqnarray*}
By the log canonical bundle formula and then by (i), we have
\begin{eqnarray*}
&&H^0(Y', \llcorner ib(K_{Y'}+\Xi)+i B^{\Xi}_-\lrcorner)\hskip 2.7in\ \\
&&\hskip 0.5in \simeq 
H^0(Y',\llcorner i(h')^*(bK_{S'}+L^{ss}_{(Y',\Xi)/S'}
+\sum s_P^{\Xi} P)+i B^{\Xi}_+\lrcorner)\\
&&\hskip 0.5in \simeq H^0(S', \llcorner
ibK_{S'}+iL_{Y'/S'}^{\log,ss}+
\sum i \, s_P^{\Xi} P\lrcorner).
\end{eqnarray*}
Thus (iii) is settled.
The property (iv) will be settled by Theorem \ref{nef} below, 
and (v) at the end of this section.     
\end{proof}

The following proposition is \cite[Proposition 4.6]{fm}. 
For the proof, see \cite{fm}. 

\begin{prop5}\label{4.4}
Under the notation and the assumptions of {\em{Theorem \ref{mainthm}}}, 
$L^{ss}_{(Y',\Xi)/S'}$ does not depend on the choice of
$\Xi$.
In particular, $L^{ss}_{(Y',\Theta')/S'}=L^{ss}_{(Y',\Theta'_+)/S'}
=L^{ss}_{(Y',\Xi)/S'}$. 
\end{prop5}

\begin{rem5}
The log canonical bundle formula:~Theorem \ref{mainthm}, 
coincides with \cite[Thoeorem 4.5]{fm} if 
$(Y',\Xi)$ is sub klt. 
So, \cite[Proposition 4.7]{fm} holds without any changes. 
\end{rem5}

The next theorem is Theorem \ref{mainthm} (iv). 

\begin{thm5}\label{nef}
The log-semistable part $L^{\log,ss}_{Y'/S'}$ is nef 
when it is restricted to a complete subvariety of $S'$, 
that is, $L^{\log,ss}_{Y'/S'}\cdot C\geq 0$ for 
every projective curve $C$ on $S'$. 
\end{thm5}
\begin{proof}
By the definition of $B^\Delta$, ${f}_{*}{\mathcal O}_{X} (\ulcorner 
(1/b)B^{\Delta}_+\urcorner)\simeq 
{\mathcal O}_{S} $ holds. 
This implies that ${h'}_*\mathcal O_{Y'}({-G'}^{<1})$ satisfies the 
condition (3) in Theorem \ref{kpt}. 
See also Claim (A) in the proof of Theorem \ref{mt}. 
By the construction \ref{(ss4.1)}, 
$G'$ 
and $\Sigma'$ satisfy the conditions (1) and (2) 
in Theorem \ref{kpt}. 
We note that $\Supp({h'}^{-1}(\Sigma')\cup G')$ is a simple normal crossing 
divisor on $Y'$. 
So, we can apply Theorem \ref{kpt} to 
$h':(Y',\Theta'-(1/b)B^{\Theta'})=(Y',G')\longrightarrow S'$. 
Thus, we obtain that $L^{\log,ss}_{Y'/S'}\cdot 
C\geq 0$ for every projective curve $C$ on $S'$. 
\end{proof}

We recall the following lemma to prove Theorem \ref{mainthm}.(v).

\begin{lem5}[{\cite[Lemma 4.12]{fm}}]\label{(4.10)}
Under the notation and the assumptions of {\em{Theorem 
\ref{mainthm} (v)}}, 
assume that $S'$ is a curve.
Then the following holds.
$$b (K_{Y'/S'}+\Xi +({(h')}^{-1}\Sigma')_{red})
\succ (h')^*(L^{\log,ss}_{Y'/S'} +b\Sigma'). 
$$
\end{lem5}

Finally, we give the proof of 
Theorem \ref{mainthm} (v) for 
the readers' convenience. 
It is essentially the same as \cite[Proof of 4.5 (i)]{fm}. 

\begin{proof}[Proof of {\em{Theorem \ref{mainthm} (v)}}]
Replacing $S'$ with a general hyperplane-section $H$ and $Y'$ 
by $(h')^*(H)$, we can immediately reduce to the case where 
$S'$ is a curve. 
For simplicity, $\Xi$ in $B^{\Xi}, s_P^{\Xi}$ and $t_P^{\Xi}$
will be suppressed during the proof.
We note that $B$ is effective.

By the hypothesis, the vertical part $D$ of the Weil divisor
$$bN (K_{Y'/S'}+\Xi^v)-(h')^*NL_{Y'/S'}^{\log,ss}=
N\sum_P s_P (h')^*P + N B-bN\Xi^h$$
is a Weil divisor.
We note that
$$
D=N\sum_P s_P (h')^*P + N B^v
=bN(K_{Y'/S'}+\Xi)-(h')^*NL_{Y'/S'}^{\log,ss}-N B^h.
$$
By Lemma \ref{(4.10)}, we have
$$
bN(K_{Y'/S'}+\Xi) 
-(h')^*NL^{\log,ss}_{Y'/S'}+bN\, ({(h')}^{-1}\Sigma')_{red}
\succ (h')^*bN\Sigma'. 
$$
Whence
\begin{equation}\label{get.ans}
D+NB^h+bN\, ((h')^{-1}\Sigma')_{red} \succ (h')^*bN\Sigma'.
\end{equation}

Let $D_P$ and $B_P^v$ be the parts of $D$ and $B^v$ lying over
$P$. 
Let $(h')^*P=\sum_k a_k F_k$ be the irreducible decomposition.
Then $D_P-NB_P^v=Ns_P (h')^*P$ and
$\Supp(D_P-Ns_P(h')^*P)\not\supset F_c$ for some $c$
by the definition of $B_P^v$. 
In particular $Ns_P a_c \in \mathbb Z$.
Furthermore, comparing the coefficients of $F_c$
in the formula (\ref{get.ans}),
we obtain
$Ns_P a_c + bN \ge bN a_c$, that is,
$Na_c s_P \ge bN(a_c-1)$.
Since $(Y',\Xi)$ is lc, we have $t_P\geq 0$ and hence $s_P\leq b$.
Hence $u_P:=a_c$ works.
\end{proof}
      
\section{Appendix:~A remark on Section \ref{sec3} by M.~Saito}
\label{sec6}
In this section, we give a different proof to Theorems 
\ref{can-ext}, \ref{loc-free}. 
It is based on the theory of mixed Hodge Modules \cite{mhm}, 
\cite{kc}. As I explained in \ref{ap}, the following \ref{6.1} 
is \cite{e-s}. 
I made no contribution in this section. 

\begin{say5}[{\cite{e-s}}]\label{6.1} 
%%%%%%%%%%%%%%%%%%%%%%%%%%%%%%%%%
%  e-mail from Morihiko Saito   %
%%%%%%%%%%%%%%%%%%%%%%%%%%%%%%%%%
Let $X$ be a smooth complex algebraic variety, and $D$ a divisor
with normal crossings whose irreducible components $D_i$ are smooth. 
Let $U=X\setminus D$ with the inclusion $j:U\to X$.
Let $(M;F,W)$ be a bifiltered (left) ${\mathcal O}_X$-Module underlying
a mixed Hodge Module. 
Assume that $L:=M|_U$ is a locally free
${\mathcal O}_U$-Module, i.e. 
it underlies an admissible variation of
mixed Hodge structure on $U$.

By the definition of pure Hodge Modules, 
we have the strict support
decomposition $$ \hbox{Gr}_k^W(M,F)=\bigoplus_Z (M_{k,Z},F),$$
where $Z$ is either $X$ or a closed irreducible variety of $D$
(by the assumption on $M|_U$), 
and the $M_{k,Z}$ have no nontrivial
subobject or quotient object with strictly smaller support.

\begin{propa}\label{p1}  
Let $p_0=\min\{p:F_pM\ne 0\}$, 
and assume $$F_{p_0}M_{k,Z}=0\quad
\hbox{if}\,\,Z\subset D.\leqno(1)$$ 
Then we have the canonical
isomorphism $$F_{p_0}M=j_*F_{p_0}L\cap L^{>-1}\leqno(2)$$
where $L^{>a}$ is the Deligne extension of $L$ such that the
eigenvalues of the residue of the connection are contained in
$(a,a+1]$.
\end{propa}
\begin{proof} 
We first consider the case $M=j_!L$, 
where $j_!$ is defined to be the composition ${\bf D}j_*{\bf D}$. 
Here ${\bf D}$
denotes the functor assigning the dual, and $j_*$ coincides with
the usual direct image as ${\mathcal O}$-Modules. 
In this case the filtration $F$ on $M$ is given 
by $$F_pM=\sum_i F_i{\mathcal D}_X
(F_{p-i}L^{>-1}),\leqno(3)$$ (see e.g.~\cite[(3.10.8)]{mhm}), where $F$
on ${\mathcal D}_X$ is the filtration by the order of operator, and
$F_pL^{>-1}$ is given by $j_*F_pL\cap L^{>-1}$ as usual.
So the isomorphism (2) is clear.

In general we use the canonical morphism $u:j_!L\to M$, see \cite
[(4.2.11)]{mhm}. 
By the above result, it is enough to show the vanishing
of $F_{p_0}$ for $\hbox{Ker}\,u$ and $\hbox{Coker}\,u$, 
because the functor assigning $F_p$ is an exact functor 
for mixed Hodge Modules.
Furthermore the functor assigning $F_p\hbox{Gr}_k^W$ is also exact.
So we may replace $u$ with $\hbox{Gr}_k^Wu:\hbox{Gr}_k^Wj_!L\to
\hbox{Gr}_k^WM$. 
This morphism is compatible with the decomposition
by strict support, and condition (1) is also satisfied for $j_!L$
(using (3)). 
So the assertion follows from the fact that
$\hbox{Gr}_k^Wu$ induces an isomorphism between the direct factors
with strict support $X$ 
(this follows from the definition of the
Hodge filtration on pure Hodge Modules, see e.g.~\cite[(3.10.12)]{mhm}.
\end{proof}

We apply this to the direct image of ${\mathcal D}$-Modules. 
Here it is easier to use right ${\mathcal D}$-Modules 
(because it simplifies many
definitions) and we use the transformation between right and left
${\mathcal D}$-Modules, 
which is defined by assigning $\Omega_X^{\dim X}
\otimes_{{\mathcal O}_X}M$ to a left ${\mathcal D}$-Module $M$, 
where $\Omega_X^{\dim X}$ has the filtration $F$ 
such that $\hbox{Gr}_{p}^F=0$ for $p\ne -\dim X$. 
We define the Hodge filtration $F$ on the right
${\mathcal D}$-Module $\omega_X$ 
by $F_p\omega_X=\omega_X$ for $p\ge 0$ 
and $0$ otherwise. Then $(\omega_X,F)$ is pure of weight $-\dim X$
(and $(\Omega_X^{\dim X},F)$ has weight $\dim X$). 
We can verify that
$\hbox{Gr}_{-k}^W(j_*\omega_U,F)$ is the direct sum of $(\iota_I)_*
(\omega_{D_I},F)$ with $\dim D_I=k$, 
where $D_I=\cap_{i\in I}D_i$ with
the inclusion $\iota_I:D_I\to X$, see 
\cite[(3.10.8) and (3.16.12)]{mhm}.
(Here the direct image $(\iota_I)_*$ is defined by tensoring the
polynomial ring ${\bf C}[\partial_1,\dots,\partial_r]$ over ${\bf C}$
if $I=\{1,\dots,r\}$, where $\partial_i=\partial/\partial x_i$ with
$(x_1,\dots,x_n)$ a local coordinate system such that $D_i=x_i^{-1}
(0)$.) 
We also see that $F_0j_*\omega_U= \omega_X(D)$, and
$$F_0H^if_*(j_*\omega_U)=R^if_*\omega_X(D)$$ by the definition of
the direct image of filtered right ${\mathcal D}$-Modules, using the
strictness of the Hodge filtration $F$ on the direct image.

\begin{propa}\label{p2} 
Let $X,U,D,j$ be as above. 
Let $f:X\to Y$ be a proper morphism of
smooth complex algebraic varieties, and $D'$ be a divisor with normal
crossings on $Y$. 
Assume that every irreducible component of any
intersections of $D_i$ is dominant to $Y$ and smooth over $Y\setminus
D'$. 
Then condition $(1)$ with $p_0=0$ is satisfied for the direct image
of a filtered $($right$)$ 
${\mathcal D}$-Module $H^if_*(j_*\omega_U,F)$.
\end{propa}
\begin{proof}
Consider the weight spectral sequence of filtered $($right$)$
${\mathcal D}$-Modules $$E_1^{-k,i+k}=H^if_*\hbox{Gr}_k^W(j_*\omega_U,F)
\Rightarrow H^if_*(j_*\omega_U,F),$$ which underlies a spectral
sequence of mixed Hodge Modules and degenerates at $E_2$. Since
$\hbox{Gr}_k^W(j_*\omega_U,F)$ is calculated as above and the direct
image of a pure Hodge Module by a proper morphism is pure, the
assertion is reduced to the proper case, where it is well known.
(Indeed, it is reduced to the torsion-freeness using the
decomposition by strict support as above.) 
\end{proof}
\end{say5}

\begin{say5}\label{dasoku}
Finally, we add one remark for the readers' convenience. 

\begin{rem5}[Deligne's extension]\label{del} 
In the above Proposition \ref{p1} and \cite[p.513]{kc}, 
$L^{>a}$ (resp.~$L^{\geq a}$) is {\em{Deligne's extension}} 
of $L$ such that the eigenvalues of the residue of 
the connection are contained in $(a, a+1]$ (resp.~$[a, a+1)$). 
In our notation:~Koll\'ar's notation \cite[Definition 2.3]{ko2}, 
$L^{>-1}$ (resp.~$L^{\geq0}$) is called 
the {\em{upper {\em{(resp.}}~lower{\em{)}} canonical extension}} of 
$L$. 
In \cite[Lemma 1.9.1]{kashi}, $L^{>-1}$ (resp.~$L^{\geq0}$) is called 
the {\em{right {\em{(resp.}}~left{\em{)}} canonical extension}} of 
$L$. 
\end{rem5}
\end{say5}
%%%%%%%%%%%%%%%%%%%%%%%%%%%%%%%%%
\ifx\undefined\bysame
\newcommand{\bysame|{leavemode\hbox to3em{\hrulefill}\,}
\fi

\end{document}